\documentclass{article} 
\usepackage{amsmath,amsthm}     
\usepackage{graphicx}     
\usepackage{hyperref} 
\usepackage{url}
\usepackage{amsfonts} 
\usepackage{multicol}
\usepackage{tikz}
\usepackage{caption}
\usepackage{subcaption}
\usepackage{stanli}
\usepackage{float}
\usepackage{hyperref}

\theoremstyle{theorem}

\theoremstyle{definition}
\newtheorem*{definition}{Definition}

\newcommand{\sator}{
\foreach \x in {1,2,3,4,5}
\foreach \y in {1,...,5}
{
  \draw (\x,\y)  rectangle ++(1,1);
}
\draw (1,1)+(0.5,0.5) node {\huge{R}}; 
\draw (2,1)+(0.5,0.5) node {\huge{O}};
\draw (3,1)+(0.5,0.5) node {\huge{T}};
\draw (4,1)+(0.5,0.5) node {\huge{A}};
\draw (5,1)+(0.5,0.5) node {\huge{S}};
\draw (1,2)+(0.5,0.5) node {\huge{O}};
\draw (2,2)+(0.5,0.5) node {\huge{P}};
\draw (3,2)+(0.5,0.5) node {\huge{E}};
\draw (4,2)+(0.5,0.5) node {\huge{R}};
\draw (5,2)+(0.5,0.5) node {\huge{A}};
\draw (1,3)+(0.5,0.5) node {\huge{T}};
\draw (2,3)+(0.5,0.5) node {\huge{E}};
\draw (3,3)+(0.5,0.5) node {\huge{N}};
\draw (4,3)+(0.5,0.5) node {\huge{E}};
\draw (5,3)+(0.5,0.5) node {\huge{T}};
\draw (1,4)+(0.5,0.5) node {\huge{A}};
\draw (2,4)+(0.5,0.5) node {\huge{R}};
\draw (3,4)+(0.5,0.5) node {\huge{E}};
\draw (4,4)+(0.5,0.5) node {\huge{P}};
\draw (5,4)+(0.5,0.5) node {\huge{O}};
\draw (1,5)+(0.5,0.5) node {\huge{S}};
\draw (2,5)+(0.5,0.5) node {\huge{A}};
\draw (3,5)+(0.5,0.5) node {\huge{T}};
\draw (4,5)+(0.5,0.5) node {\huge{O}};
\draw (5,5)+(0.5,0.5) node {\huge{R}};
}

\allowdisplaybreaks

\makeatletter
\@addtoreset{footnote}{page}
\makeatother

\begin{document}
\usetikzlibrary {arrows.meta}  

\title{Some mathematical and geometrical interpretations of the Sator Square}

\author{Paul Dario TOASA CAIZA\\               
\small Karlsruher Institut für Technologie. Versuchsanstalt für Stahl, Holz
und Steine\\    
\small Otto-Ammann-Platz 1. 76131 Karlsruhe. Germany\\                
\small paul.toasa@kit.edu}                      

\maketitle

\noindent
In 1738, the King of Naples and future King of Spain, Carlos III, commissioned
the Spanish military engineer Roque Joaqu\'in de Alcubierre to begin the
excavations of the ruins of the ancient Roman city of Pompeii and its
surroundings, buried by the terrible explosion of Vesuvius in AD 79. 
Since that time, archaeologists have brought to light wonderful treasures
found in the among ruins.
Among them, the Sator Square is one of the most peculiar, apparently simple
but mysterious. 
Supernatural and medicinal powers have been attributed to this object and its
use was widespread during the Middle Age. 
Studies to explain its origin and meaning have been varied.
There are theories that relate it to religion, the occult, medicine and music.
However, no explanation has been convincing beyond pseudo-scientific sensationalism.
In this study, the author intends to eliminate the mystical character of the
Sator Square and suggests considering it as a simple palindrome or a game of
words with certain symmetrical properties.
However, these properties are not exclusive to the Sator Suare but are present in
various mathematical and geometric objects. 

\section{Introduction}

I was very intrigued by an article I read about the enigmatic
\href{https://www.bbc.com/mundo/noticias-64188793}{Sator Square} on the BBC
Spanish service.  
And the intrigue was because all the mentioned explanations were spiritual and
pseudo-scientific. 
Despite the fact that the object presented certain symmetries, no geometric or 
mathematical study was mentioned in this regard. 
So I searched for mathematical studies on this object and only found countless
sensational and absurd works. 
So I asked myself: Is it possible to consider the Sator Square as a mathematical
object and study it? \\

\begin{figure}[h!]
  \centering
\begin{tikzpicture}
\sator
\end{tikzpicture}
\caption{Sator square}
\label{fig:sator}
\end{figure}

Among the many wonderful historical artifacts found in the ruins of Pompeii
caused by the eruption of Vesuvius in A.D. 79, Sator Square is one of the most
unique and unusual, see Figure \ref{fig:sator}.  
Particularly, during the Middle Age, the square was used as talisman or amulet,
against rabies and other diseases and complaints, such as fever, intestinal
gout and rose, jaundice and consumption, cramps, headache and toothache, birth
complications and worms.
It was also regarded as a remedy against thieves, witchcraft and fire, as a
spell of exorcism, blessing of bullets and hunting aid, even as a 
love magic \cite{bader:1987}.

An immense amount of works have been written in order to explain its origin,
meaning and, in the case of being a riddle, its solution.
All these works suggest a relationship with Christianity, Judaism, Satanism,
philology, numerology, esotericism, magic, mysticism and even medicine
\cite{seligmann:1914}, \cite{darmstaedter:1932}, \cite{fishwick:1964}.
However, none of these works provides a complete and convincing explanation
about the origin and meaning of the square.\\

The magic of this square rest surely on the perfect symmetry of its component
letters which yield the same combinations in four different directions
\cite{fishwick:1964}.
It has 5 Latin words, which can be read from left to right and from top to bottom,
or from left to right and bottom to top and their ``message'' is the same, see
Figure \ref{fig:palindrome}

\begin{figure}[h!]
  \centering
  \Large{SATOR AREPO TENET OPERA ROTAS}
  \caption{Five words palindrome}
  \label{fig:palindrome}
\end{figure}

In \cite{fishwick:1964}, the author analyses how difficult to build a rebus
with twenty five letters in five Latin words is.
Moreover, it is suggested that the words were written by Latin-speaking Jews.\\
Other similar squares have been written by considering other words and other
languages like Greek and Hebrew see \cite{seligmann:1914}.\\
Even in the pop culture the Sator Square has inspired the creation of Tenet,
which is a 2020 science fiction action thriller film directed and written by
Christopher Nolan.
In this film the Protagonist master the art of "time inversion" as a way of
countering the threat that is to come.\\
In \cite{griffiths:1971} the origin of AREPO is discussed and it is proposed
that the square has Egyptian origin.\\
In \cite{ferguson:1974} the work of Moeller \cite{moeller:1973} is reviewed.
Moeller proposes that Sator is Saturn and assures "I now present these
findings in the firm conviction that I have come upon the correct solution``.
This declaration received a lot of critisism, saying that this work went in
the wrong direction.
Others called it provocative but no persuasive.\\
In fact, about half the text deals with numerology and it also suggests that
probably the square has Christian origin, fact that Hemer supports in
\cite{hemer:1978}.
This is based on the fact that the phrase ``Pater Noster'' and additionally the
letters A and O, standing for Alpha and Omega, the first and the last were obtained by
rearrangin the twenty five letters of the square.\\
In \cite{hofmann:1976}, the work of Moeller is compared to a sandcastle.\\

Even in the classical music the influence of the square can be noted.
The Austrian composer Anton Friedrich Wilhelm von Webern\footnote{3.12.1883 –
  15.09.1945} considered the square strucuture in its composition ``Concerto
for Nine Instruments, Op. 24''\footnote{Konzert f\"ur neun Instrumente}, written
in 1934 \cite{neuwirth:1980}. \\

An extensive study about several interpretations of the square, with a
abundant bibliography is available in \cite{sheldon:2003}. 

\section{Reading simmetry}
The most discussed chracteristic of the Sator square is related with the
concept of palyndromes.
The Sator square has five words, located in a particular way, which makes
possible to obtain several types of reading simmetry.

\subsection{Clasic palindrome}
\begin{figure}[h!]
  \centering

\begin{tikzpicture}
\sator
\draw[red] (3.9,6.8) node {$\mathbf{S_y}$};
\draw[red,ultra thick,{Stealth[length=3mm]}-{Stealth[length=3mm]}]
(0,3.5)--(7,3.5);
\draw[red] (6.8, 3.9) node {$\mathbf{S_x}$};
\draw[red,ultra thick,{Stealth[length=3mm]}-{Stealth[length=3mm]}] (3.5,7)--(3.5,0);
\end{tikzpicture}
\caption{The palindrome TENET of the Sator square}
\label{fig:tenet}
\end{figure}
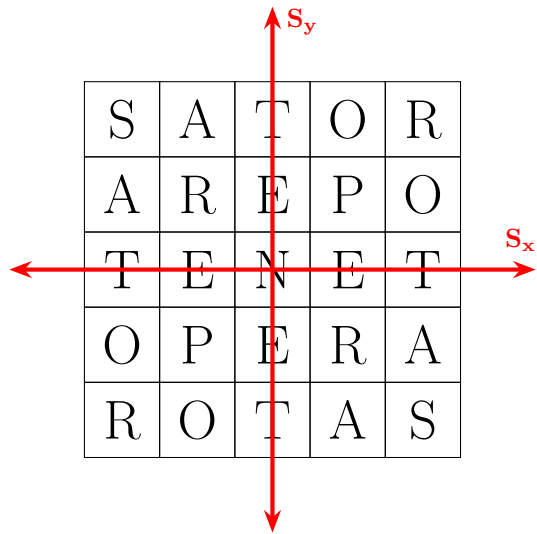

The word TENET is a palyndrome, and it is written in the  middle column and
row of the square.
If two perpendicular axis are drawn trought the center of the square,
the palindrome TENET coincides with these axes, see Figure \ref{fig:tenet}.

Moreover, by writing consecutively the five words, a longer palindrome is
obtained as well, see Figure \ref{fig:palindrome}.

\subsection{Palindrome in both directions}
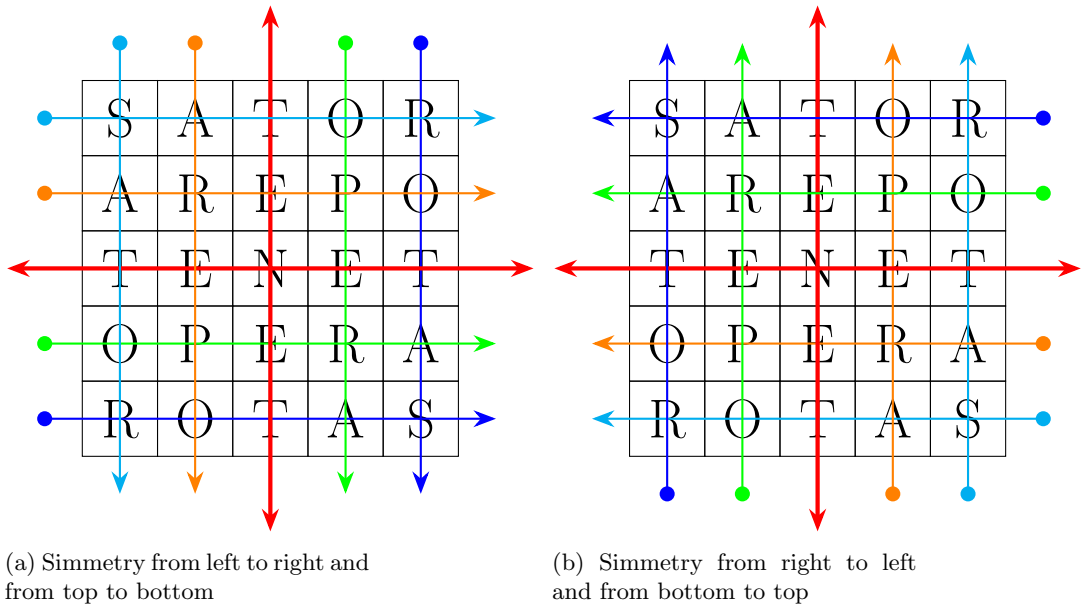
\begin{figure}[h!]
\centering
\begin{subfigure}[b]{0.4\textwidth}
\begin{tikzpicture}
\sator
\fill[blue] (0.5,1.5) circle(1mm);
\draw[blue,thick,-{Stealth[length=3mm]}] (0.5,1.5)--(6.5,1.5);
\fill[blue] (5.5,6.5) circle(1mm);
\draw[blue,thick,-{Stealth[length=3mm]}] (5.5,6.5)--(5.5,0.5);
\fill[green] (0.5,2.5) circle(1mm);
\draw[green,thick,-{Stealth[length=3mm]}] (0.5,2.5)--(6.5,2.5);
\fill[green] (4.5,6.5) circle(1mm);
\draw[green,thick,-{Stealth[length=3mm]}] (4.5,6.5)--(4.5,0.5);
\draw[red,ultra thick,{Stealth[length=3mm]}-{Stealth[length=3mm]}] (0,3.5)--(7,3.5);
\draw[red,ultra thick,{Stealth[length=3mm]}-{Stealth[length=3mm]}] (3.5,7)--(3.5,0);
\fill[orange] (0.5,4.5) circle(1mm);
\draw[orange,thick,-{Stealth[length=3mm]}] (0.5,4.5)--(6.5,4.5);
\fill[orange] (2.5,6.5) circle(1mm);
\draw[orange,thick,-{Stealth[length=3mm]}] (2.5,6.5)--(2.5,0.5);
\fill[cyan] (0.5,5.5) circle(1mm);
\draw[cyan,thick,-{Stealth[length=3mm]}] (0.5,5.5)--(6.5,5.5);
\fill[cyan] (1.5,6.5) circle(1mm);
\draw[cyan,thick,-{Stealth[length=3mm]}] (1.5,6.5)--(1.5,0.5);
\end{tikzpicture}
\caption{Simmetry from left to right and from top to bottom}
\label{sfig:lr}
\end{subfigure}
\hfill
\begin{subfigure}[b]{0.4\textwidth}
\begin{tikzpicture}
\sator
\fill[cyan] (6.5,1.5) circle(1mm);
\draw[cyan,thick,{Stealth[length=3mm]}-] (0.5,1.5)--(6.5,1.5);
\fill[cyan] (5.5,0.5) circle(1mm);
\draw[cyan,thick,{Stealth[length=3mm]}-] (5.5,6.5)--(5.5,0.5);
\fill[orange] (6.5,2.5) circle(1mm);
\draw[orange,thick,{Stealth[length=3mm]}-] (0.5,2.5)--(6.5,2.5);
\fill[orange] (4.5,0.5) circle(1mm);
\draw[orange,thick,{Stealth[length=3mm]}-] (4.5,6.5)--(4.5,0.5);
\draw[red,ultra thick,{Stealth[length=3mm]}-{Stealth[length=3mm]}] (0,3.5)--(7,3.5);
\draw[red,ultra thick,{Stealth[length=3mm]}-{Stealth[length=3mm]}] (3.5,7)--(3.5,0);
\fill[green] (6.5,4.5) circle(1mm);
\draw[green,thick,{Stealth[length=3mm]}-] (0.5,4.5)--(6.5,4.5);
\fill[green] (2.5,0.5) circle(1mm);
\draw[green,thick,{Stealth[length=3mm]}-] (2.5,6.5)--(2.5,0.5);
\fill[blue] (6.5,5.5) circle(1mm);
\draw[blue,thick,{Stealth[length=3mm]}-] (0.5,5.5)--(6.5,5.5);
\fill[blue] (1.5,0.5) circle(1mm);
\draw[blue,thick,{Stealth[length=3mm]}-] (1.5,6.5)--(1.5,0.5);
\end{tikzpicture}
\caption{Simmetry from right to left and from bottom to top}
\label{sfig:rl}
\end{subfigure}
\caption{Reading simmetry of the Sator square in two ways}
\label{fig:lrtb}
\end{figure}

The four words SATOR, AREPO, OPERA and ROTAS behave like a palindrome. 
They can be observed and read in the traditional western way
from left to right and from top to bottom or from right to left and from
bottom to top.
Depending on the reading direction these words are located in the ith-row and
ith-column.
For example, OPERA is located in the 4th row and the 4th column from Figure
\ref{sfig:lr} and in the 2nd row and 2nd column  from Figure \ref{sfig:rl},
see the corresponding green arrows.

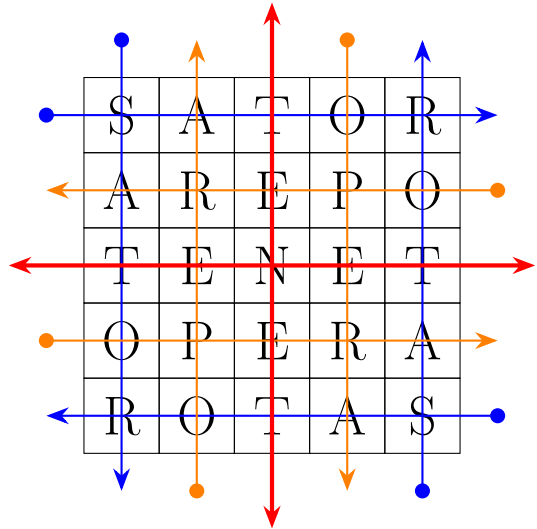
\begin{figure}[h!]
\centering
\begin{tikzpicture}
\sator
\fill[blue] (0.5,5.5) circle(1mm);
\draw[blue,thick,-{Stealth[length=3mm]}] (0.5,5.5)--(6.5,5.5);
\fill[blue] (5.5,0.5) circle(1mm);
\draw[blue,thick,{Stealth[length=3mm]}-] (5.5,6.5)--(5.5,0.5);
\fill[blue] (6.5,1.5) circle(1mm);
\draw[blue,thick,-{Stealth[length=3mm]}] (6.5,1.5)--(0.5,1.5);
\fill[blue] (1.5,6.5) circle(1mm);
\draw[blue,thick,-{Stealth[length=3mm]}] (1.5,6.5)--(1.5,0.5);
\fill[orange] (0.5,2.5) circle(1mm);
\draw[orange,thick,-{Stealth[length=3mm]}] (0.5,2.5)--(6.5,2.5);
\fill[orange] (4.5,6.5) circle(1mm);
\draw[orange,thick,-{Stealth[length=3mm]}] (4.5,6.5)--(4.5,0.5);
\fill[orange] (6.5,4.5) circle(1mm);
\draw[orange,thick,{Stealth[length=3mm]}-] (0.5,4.5)--(6.5,4.5);
\fill[orange] (2.5,0.5) circle(1mm);
\draw[orange,thick,{Stealth[length=3mm]}-] (2.5,6.5)--(2.5,0.5);
\draw[red,ultra thick,{Stealth[length=3mm]}-{Stealth[length=3mm]}] (0,3.5)--(7,3.5);
\draw[red,ultra thick,{Stealth[length=3mm]}-{Stealth[length=3mm]}] (3.5,7)--(3.5,0);
\end{tikzpicture}
\caption{Reading simmetry of the Sator square in four ways}
\label{fig:4ways}
\end{figure}

If both squares from Figure \ref{fig:lrtb} are joined, the five words of the
Sator square are reduced to three SATOR, AREPO and TENET, see Figure
\ref{fig:4ways}.
Moreover, it is also possible to observe that the Sator square keeps its form
after rotating one time  around the horizontal axis and one time around the
vertical axis.
This kind of simmetry will be discussed  in the next section.

\clearpage


\section{Geometric symmetry}
The reading symmmetries observed in the previous section can be described
trought different types of geometrical symmetries. 
This description will allow to find a general way to generate Sator Squares of
higher dimensions.

\subsection{Point symmetry}
When every part of an object has a matching part, and both parts have the same 
distance but in the opposite direction from the centre point, this is known as
point symmetry. 
It’s also known as “Order 2 rotational symmetry.”\\
Because the “Origin” is the primary point around which the shape is
symmetrical, point symmetry is also known as origin symmetry.\\
When a figure is rotated $180^o$, and it keeps its original shape, it is said
that it has point symmetry.\\
Thus, the palindrome TENET in the Sator Square has a point symmetry located at the
letter N, see Figure \ref{fig:tenet}.
In other words, the axes $S_x$ and $S_y$ have a point symmetry with
respect to the center of the square.

\subsection{Reflective or mirror symmetry}
\begin{figure}[h!]
  \centering
\begin{tikzpicture}
\sator
\draw[red] (6.7,6.7) node {$\mathbf{S_1}$};
\draw[red,ultra thick,dashed,{Stealth[length=3mm]}-{Stealth[length=3mm]}](0.5,0.5)--(6.5,6.5);
\draw[blue] (0.3,6.7) node {$\mathbf{S_2}$};
\draw[blue,ultra thick,dashed,{Stealth[length=3mm]}-{Stealth[length=3mm]}](0.5,6.5)--(6.5,0.5); 
\end{tikzpicture}
\caption{Lines of simmetry of the Sator square}
\label{fig:S1S2}
\end{figure}
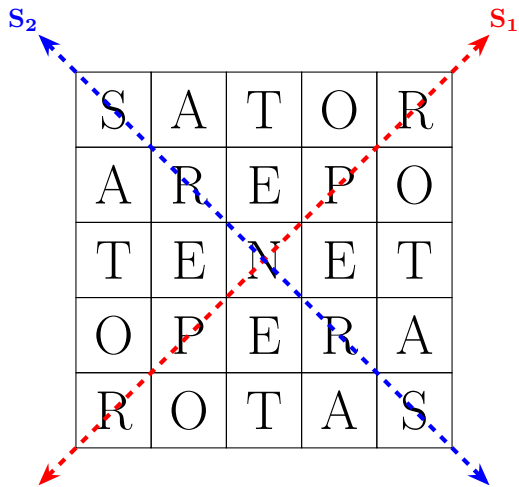

Reflective symmetry, also called mirror symmetry, is a type of symmetry where
one half of the object reflects the other half of the object.
In other words, in this type of symmetry one part of the image or object
represents the mirror image of another part of the image.
Reflexive Symmetry is described by a line of symmetry.\\
In the Sator Square two lines of symmetry $S_1$ and $S_2$ are present, see
Figure \ref{fig:S1S2}.

\begin{figure}[h!]
  \centering
  \begin{subfigure}[b]{0.45\textwidth}
\begin{tikzpicture}
  \sator
\draw[blue] (0.3,6.7) node {$\mathbf{S_2}$};  
\draw[red,ultra thick,dashed,{Stealth[length=3mm]}-{Stealth[length=3mm]}] (0.5,3.5)--(3.5,6.5);
\draw[red,ultra thick,dashed,{Stealth[length=3mm]}-{Stealth[length=3mm]}] (0.5,2.5)--(4.5,6.5);
\draw[red,ultra thick,dashed,{Stealth[length=3mm]}-{Stealth[length=3mm]}] (0.5,1.5)--(5.5,6.5);
\draw[red,ultra thick,{Stealth[length=3mm]}-{Stealth[length=3mm]}] (0.5,0.5)--(6.5,6.5);
\draw[blue,ultra thick,{Stealth[length=3mm]}-{Stealth[length=3mm]}] (0.5,6.5)--(6.5,0.5);
\end{tikzpicture}
\caption{Reflexions  relative to symmetry line S2}
\end{subfigure}
\hfill
\begin{subfigure}[b]{0.45\textwidth}
\begin{tikzpicture}  
  \sator
\draw[red] (6.7,6.7) node {$\mathbf{S1}$};  
\draw[red,ultra thick,{Stealth[length=3mm]}-{Stealth[length=3mm]}] (0.5,0.5)--(6.5,6.5);
\draw[blue,ultra thick,dashed,{Stealth[length=3mm]}-{Stealth[length=3mm]}]
(6.5,3.5)--(3.5,6.5);
\draw[blue,ultra thick,dashed,{Stealth[length=3mm]}-{Stealth[length=3mm]}] (6.5,2.5)--(2.5,6.5);
\draw[blue,ultra thick,dashed,{Stealth[length=3mm]}-{Stealth[length=3mm]}] (6.5,1.5)--(1.5,6.5);
\draw[blue,ultra thick,{Stealth[length=3mm]}-{Stealth[length=3mm]}] (0.5,6.5)--(6.5,0.5);
\end{tikzpicture}
\caption{Reflexions  relative to symmetry line S1}
\end{subfigure}
\caption{Simmetry lines and their corresponding reflection in the Sator
  Square}
\label{fig:S1S2Reflexion}
\end{figure}
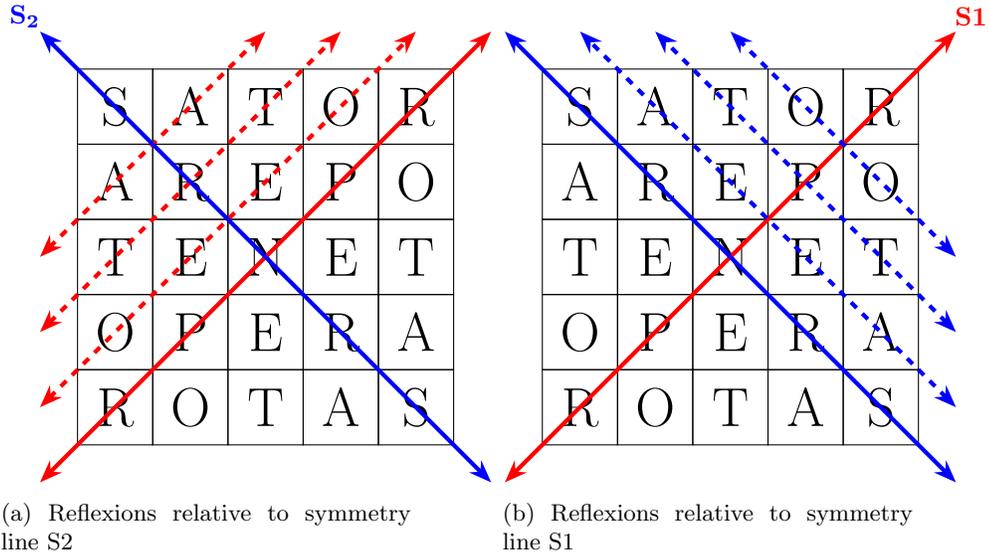

If a mirror is located along the lines of symmetry $S_1$ and $S_2$, the
letters in the Sator Square will be reflected as it is shown with dotted lines
in Figure \ref{fig:S1S2Reflexion}

\subsection{Rotational simmetry in 3D}
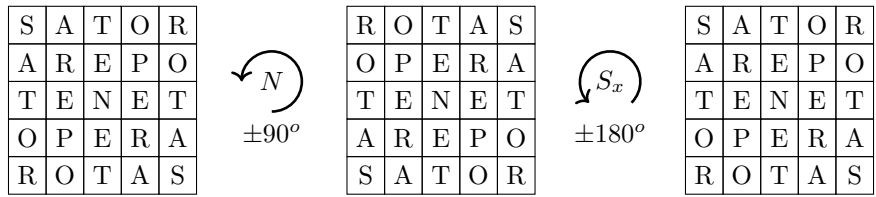
\begin{figure}[h!]
  \centering
\begin{tikzpicture}
\foreach \x in {1,2,3,4,5}
\foreach \y in {1,2,3,4,5}
{
  \draw (\x/2,\y/2)  rectangle ++(0.5,0.5);
  \draw (\x/2+4.5,\y/2) rectangle ++(0.5,0.5);
  \draw (\x/2+9,\y/2) rectangle ++(0.5,0.5);
}
\draw (0.5,0.5)+(0.25,0.25) node {R}; 
\draw (1,0.5)+(0.25,0.25) node {O};
\draw (1.5,0.5)+(0.25,0.25) node {T};
\draw (2,0.5)+(0.25,0.25) node {A};
\draw (2.5,0.5)+(0.25,0.25) node {S};
\draw (0.5,1)+(0.25,0.25) node {O};
\draw (1,1)+(0.25,0.25) node {P};
\draw (1.5,1)+(0.25,0.25) node {E};
\draw (2,1)+(0.25,0.25) node {R};
\draw (2.5,1)+(0.25,0.25) node {A};
\draw (0.5,1.5)+(0.25,0.25) node {T};
\draw (1,1.5)  +(0.25,0.25) node {E};
\draw (1.5,1.5)+(0.25,0.25) node {N};
\draw (2,1.5)  +(0.25,0.25) node {E};
\draw (2.5,1.5)+(0.25,0.25) node {T};
\draw (0.5,2)+(0.25,0.25) node {A};
\draw (1,2)  +(0.25,0.25) node {R};
\draw (1.5,2)+(0.25,0.25) node {E};
\draw (2,2)  +(0.25,0.25) node {P};
\draw (2.5,2)+(0.25,0.25) node {O};
\draw (0.5,2.5)+(0.25,0.25) node {S};
\draw (1,2.5)  +(0.25,0.25) node {A};
\draw (1.5,2.5)+(0.25,0.25) node {T};
\draw (2,2.5)  +(0.25,0.25) node {O};
\draw (2.5,2.5)+(0.25,0.25) node {R};
\point {a}{4}{2}; 
\load {3}{a}[-90];
\draw (4,1.3) node {$\pm90^o$};
\draw (4,2) node {$N$};
\draw (5,0.5)+(0.25,0.25) node {S}; 
\draw (5.5,0.5)+(0.25,0.25) node {A};
\draw (6,0.5)+(0.25,0.25) node {T};
\draw (6.5,0.5)+(0.25,0.25) node {O};
\draw (7,0.5)+(0.25,0.25) node {R};
\draw (5,1)+(0.25,0.25) node {A};
\draw (5.5,1)+(0.25,0.25) node {R};
\draw (6,1)+(0.25,0.25) node {E};
\draw (6.5,1)+(0.25,0.25) node {P};
\draw (7,1)+(0.25,0.25) node {O};
\draw (5,1.5)+(0.25,0.25) node {T};
\draw (5.5,1.5)  +(0.25,0.25) node {E};
\draw (6,1.5)+(0.25,0.25) node {N};
\draw (6.5,1.5)  +(0.25,0.25) node {E};
\draw (7,1.5)+(0.25,0.25) node {T};
\draw (5,2)+(0.25,0.25) node {O};
\draw (5.5,2)  +(0.25,0.25) node {P};
\draw (6,2)+(0.25,0.25) node {E};
\draw (6.5,2)  +(0.25,0.25) node {R};
\draw (7,2)+(0.25,0.25) node {A};
\draw (5,2.5)+(0.25,0.25) node {R};
\draw (5.5,2.5)  +(0.25,0.25) node {O};
\draw (6,2.5)+(0.25,0.25) node {T};
\draw (6.5,2.5)  +(0.25,0.25) node {A};
\draw (7,2.5)+(0.25,0.25) node {S};
--------------- Arrow
\point {b}{8.5}{2}; 
\load {3}{b}[-45];
\draw (8.5,1.3) node {$\pm180^o$};
\draw (8.5,2) node {$S_x$};
\draw (9.5,0.5)+(0.25,0.25) node {R}; 
\draw (10,0.5)+(0.25,0.25) node {O};
\draw (10.5,0.5)+(0.25,0.25) node {T};
\draw (11,0.5)+(0.25,0.25) node {A};
\draw (11.5,0.5)+(0.25,0.25) node {S};
\draw (9.5,1)+(0.25,0.25) node {O};
\draw (10,1)+(0.25,0.25) node {P};
\draw (10.5,1)+(0.25,0.25) node {E};
\draw (11,1)+(0.25,0.25) node {R};
\draw (11.5,1)+(0.25,0.25) node {A};
\draw (9.5,1.5)+(0.25,0.25) node {T};
\draw (10,1.5)  +(0.25,0.25) node {E};
\draw (10.5,1.5)+(0.25,0.25) node {N};
\draw (11,1.5)  +(0.25,0.25) node {E};
\draw (11.5,1.5)+(0.25,0.25) node {T};
\draw (9.5,2)+(0.25,0.25) node {A};
\draw (10,2)  +(0.25,0.25) node {R};
\draw (10.5,2)+(0.25,0.25) node {E};
\draw (11,2)  +(0.25,0.25) node {P};
\draw (11.5,2)+(0.25,0.25) node {O};
\draw (9.5,2.5)+(0.25,0.25) node {S};
\draw (10,2.5)  +(0.25,0.25) node {A};
\draw (10.5,2.5)+(0.25,0.25) node {T};
\draw (11,2.5)  +(0.25,0.25) node {O};
\draw (11.5,2.5)+(0.25,0.25) node {R};
\end{tikzpicture}
\caption{Rotational simmetry in 3D of the Sator Square}
\label{fig:S3D}
\end{figure}

Besides the 2D symmetries observed in the Sator Square, a symmetry in 3D is
presented as well.
If the Sator Squere is rotated $\pm90^o$, and afterwards it is rotated $\pm180^o$ by
considering as rotation axis the symmetry line $S_x$ or $S_y$, the Sator
Square keeps its original form, see Figure \ref{fig:S3D}.
Doing these operation in a reverse order, the original form is kept as well.

\section{Satorian functions and squares}

Considering as reference the reading and geometrical simmetries observed in the
Sator square, it is possible to find a general discrete way to create Sator
squares of higher dimension.

However, before proposing a general definition and for
didactic purposes, let us consider the following function, which satisfies the
point simmetry along $S_x$ and $S_y$ corresponding to the center of the square
and the lines of symmetry $S_1$ and $S_2$. 

\begin{equation}
  f(x,y)=\left\{
    \begin{array}{rl}
      \sqrt{x^2+y^2} & \text{if } xy\geq 0,\\\\
      -\sqrt{x^2+y^2} & \text{else},
    \end{array} \right.
  \label{eq:norm}
 \end{equation}
 where $x,y\in \mathbb{Z}$.

\begin{figure}[H]
\centering
\begin{subfigure}[b]{0.44\textwidth}
\begin{tikzpicture}
\foreach \x in {1,2,3,4,5}
\foreach \y in {1,...,5}
{
  \draw (\x,\y)  rectangle ++(1,1);
}
\draw (1,1)+(0.5,0.5) node {$\cdot$}; 
\draw (2,1)+(0.5,0.5) node {$\cdot$};
\draw (3,1)+(0.5,0.5) node {-2};
\draw (4,1)+(0.5,0.5) node {$\cdot$};
\draw (5,1)+(0.5,0.5) node {$\cdot$};
\draw (1,2)+(0.5,0.5) node {$\cdot$};
\draw (2,2)+(0.5,0.5) node {$\cdot$};
\draw (3,2)+(0.5,0.5) node {-1};
\draw (4,2)+(0.5,0.5) node {$\cdot$};
\draw (5,2)+(0.5,0.5) node {$\cdot$};
\draw (1,3)+(0.5,0.5) node {-2};
\draw (2,3)+(0.5,0.5) node {-1};
\draw (3,3)+(0.5,0.5) node {0};
\draw (4,3)+(0.5,0.5) node {1};
\draw (5,3)+(0.5,0.5) node {2};
\draw (1,4)+(0.5,0.5) node {$\cdot$};
\draw (2,4)+(0.5,0.5) node {$\cdot$};
\draw (3,4)+(0.5,0.5) node {1};
\draw (4,4)+(0.5,0.5) node {$\cdot$};
\draw (5,4)+(0.5,0.5) node {$\cdot$};
\draw (1,5)+(0.5,0.5) node {$\cdot$};
\draw (2,5)+(0.5,0.5) node {$\cdot$};
\draw (3,5)+(0.5,0.5) node {2};
\draw (4,5)+(0.5,0.5) node {$\cdot$};
\draw (5,5)+(0.5,0.5) node {$\cdot$};
\draw[red,thick,{Stealth[length=3mm]}-{Stealth[length=3mm]}] (0,3.5)--(7,3.5);
\draw[red, thick,{Stealth[length=3mm]}-{Stealth[length=3mm]}] (3.5,7)--(3.5,0);
\end{tikzpicture}
\caption{$(x,y)$}
\label{sfig:xy}
\end{subfigure}
\hfill
\begin{subfigure}[b]{0.44\textwidth}
\begin{tikzpicture}
\foreach \x in {1,2,3,4,5}
\foreach \y in {1,...,5}
{
  \draw (\x,\y)  rectangle ++(1,1);
}
\draw (1,1)+(0.5,0.5) node {$\sqrt{8}$}; 
\draw (2,1)+(0.5,0.5) node {$\sqrt{5}$};
\draw (3,1)+(0.5,0.5) node {2};
\draw (4,1)+(0.5,0.5) node {$-\sqrt{5}$};
\draw (5,1)+(0.5,0.5) node {$-\sqrt{8}$};
\draw (1,2)+(0.5,0.5) node {$\sqrt{5}$};
\draw (2,2)+(0.5,0.5) node {$\sqrt{2}$};
\draw (3,2)+(0.5,0.5) node {$1$};
\draw (4,2)+(0.5,0.5) node {$-\sqrt{2}$};
\draw (5,2)+(0.5,0.5) node {$-\sqrt{5}$};
\draw (1,3)+(0.5,0.5) node {$2$};
\draw (2,3)+(0.5,0.5) node {$1$};
\draw (3,3)+(0.5,0.5) node {$0$};
\draw (4,3)+(0.5,0.5) node {$1$};
\draw (5,3)+(0.5,0.5) node {$2$};

\draw (1,4)+(0.5,0.5) node {$-\sqrt{5}$};
\draw (2,4)+(0.5,0.5) node {$-\sqrt{2}$};
\draw (3,4)+(0.5,0.5) node {$1$};
\draw (4,4)+(0.5,0.5) node {$\sqrt{2}$};
\draw (5,4)+(0.5,0.5) node {$\sqrt{5}$};
\draw (1,5)+(0.5,0.5) node {$-\sqrt{8}$};
\draw (2,5)+(0.5,0.5) node {$-\sqrt{5}$};
\draw (3,5)+(0.5,0.5) node {$2$};
\draw (4,5)+(0.5,0.5) node {$\sqrt{5}$};
\draw (5,5)+(0.5,0.5) node {$\sqrt{8}$};
\draw[red,thick,{Stealth[length=3mm]}-{Stealth[length=3mm]}] (0,3.5)--(7,3.5);
\draw[red,thick,{Stealth[length=3mm]}-{Stealth[length=3mm]}] (3.5,7)--(3.5,0);
\end{tikzpicture}
\caption{$f(x,y)$}
\label{sfig:fxy}
\end{subfigure}
\caption{Application of the Equation (\ref{eq:norm}) for $|x|,|y|\leq 2$}
\label{fig:app-n5}
\end{figure}
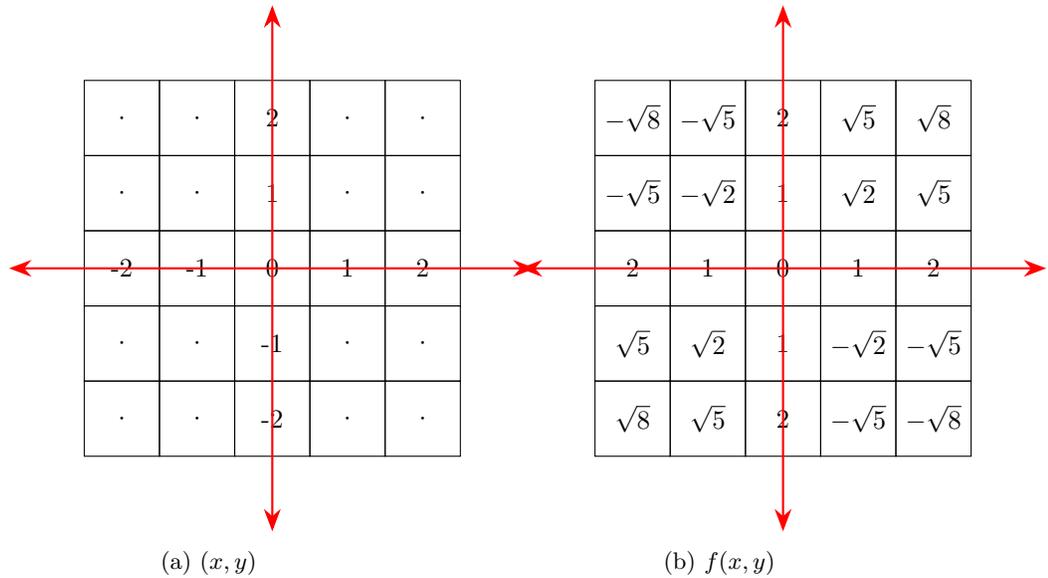

\begin{figure}[h!]
  \centering  
\begin{subfigure}[b]{0.45\textwidth}
\includegraphics[width=\textwidth]{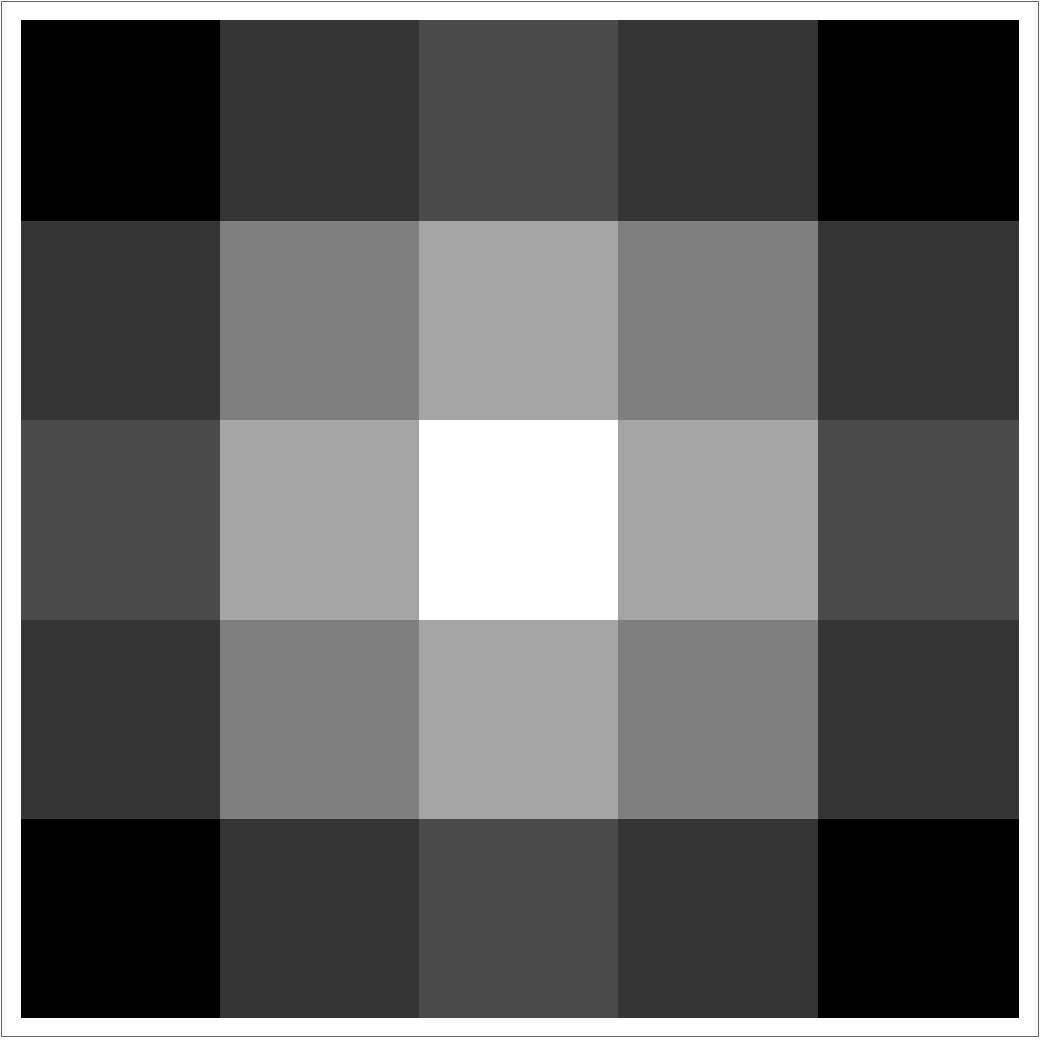}
\caption{$SS_5$ in gray scale}
\label{sfig:SS5-BW}
\end{subfigure}
\hfill
\begin{subfigure}[b]{0.45\textwidth}
  \includegraphics[width=\textwidth]{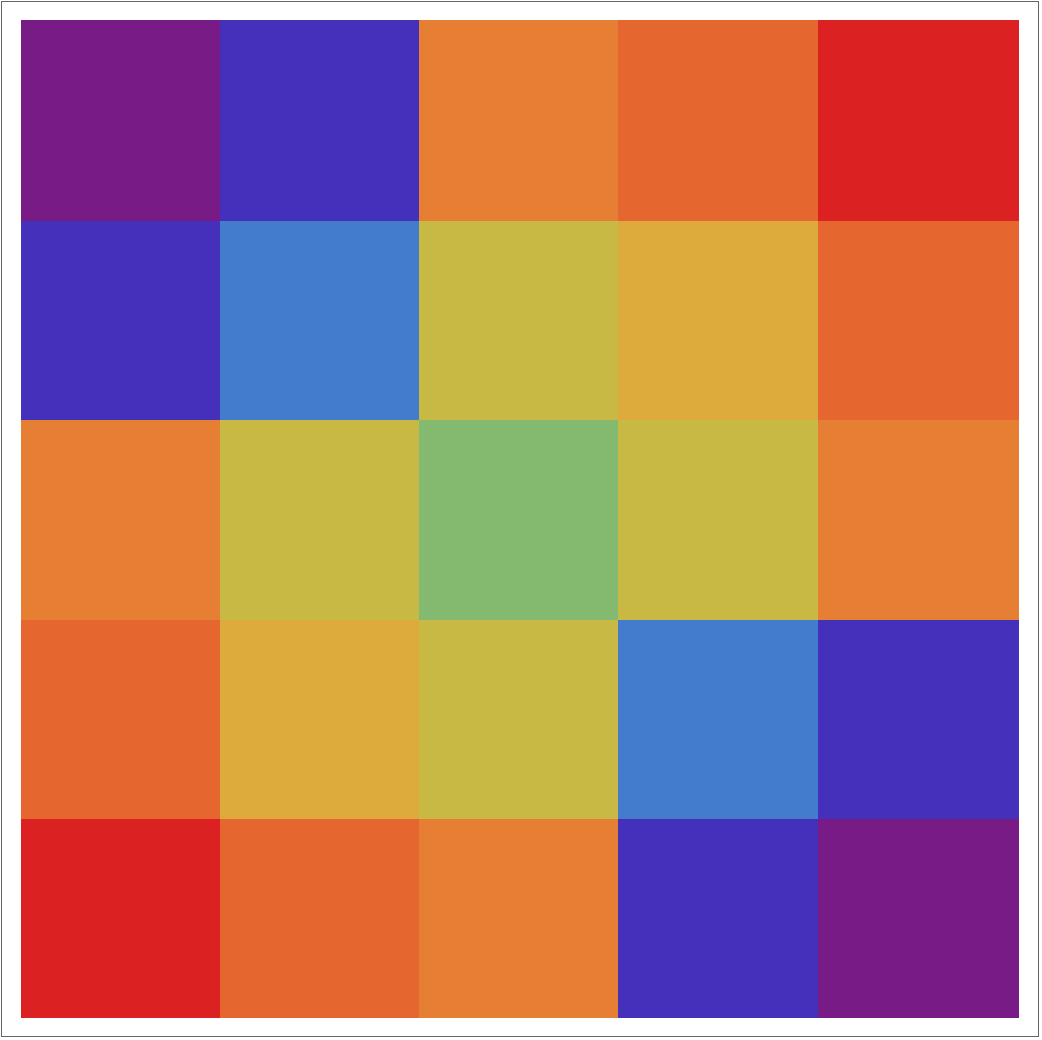}
  \caption{$SS_5$ in color scale}
  \label{sfig:SS5-C}
\end{subfigure}
\caption{Graphical representation of a Satorian square $SS_5$ in gray scale
  and color}
\label{fig:SQColor}
\end{figure}

Figure \ref{fig:app-n5} shows the application of the Equation (\ref{eq:norm})
in the set given by $|x|,|y|\leq 2$.
Consequently, a square of dimension $n=5$ similar to the Sator Square is
obtained.
In fact, the only difference is that instead of letters, numbers are used.
It is worth mentioning that in order to keep the expected simmetries, the
dimension of the generated square has to be odd.\\
By assignng colors to the values of the square in Figure \ref{sfig:fxy}, it is
possible to observe the expected simmetries in a easier way, see Figure \ref{fig:SQColor}.

The Equation (\ref{eq:norm}) is a particular case of a function, which will be
called Satorian function and it will be defined as follows.

\begin{definition}[Satorian function]
Let $(x,y)$ be a number in $\mathbb{Z}^2$, a function $f(x,y)$ is called
Satorian when $f(x,y)=f(y,x)$ and $f(x,y)=f(-y, -x)$.
\end{definition}

The squares shown in Figure \ref{fig:SQColor} are called Satorian squares and
are defined as follows.

\begin{definition}[Satorian square]
The square obtained by applying a Satorian function $f(x,y)$ on a   
a discrete, centered and symmetric set of numbers of dimension $n\times n$,
where $n$ is odd, is called Satorian square of dimension $n$ and it is noted
as $SS_n$. 
\end{definition}

It is not necessary to define a Satorian function as an explicyt equation, it
can be also given by a set of functions or by a set of discrete values.
In fact, the archeological Satorian square is given in this way.

\subsection{Archeological Satorian square}
\begin{figure}[h!]
  \centering  
\begin{subfigure}[b]{0.55\textwidth}
\begin{tikzpicture}
\foreach \x in {1,2,3,4,5}
\foreach \y in {1,...,5}
{
  \draw (\x,\y)  rectangle ++(1,1);
}
\draw (1,1)+(0.5,0.5) node {$\sqrt{8}$}; 
\draw (2,1)+(0.5,0.5) node {$\sqrt{5}$};
\draw (3,1)+(0.5,0.5) node {2};
\draw (4,1)+(0.5,0.5) node {$-\sqrt{5}$};
\draw (5,1)+(0.5,0.5) node {$-\sqrt{8}$};
\draw (1,2)+(0.5,0.5) node {$\sqrt{5}$};
\draw (2,2)+(0.5,0.5) node {$\sqrt{2}$};
\draw (3,2)+(0.5,0.5) node {$1$};
\draw[blue] (4,2)+(0.5,0.5) node {$\sqrt{8}$};
\draw (5,2)+(0.5,0.5) node {$-\sqrt{5}$};
\draw (1,3)+(0.5,0.5) node {$2$};
\draw (2,3)+(0.5,0.5) node {$1$};
\draw (3,3)+(0.5,0.5) node {$0$};
\draw (4,3)+(0.5,0.5) node {$1$};
\draw (5,3)+(0.5,0.5) node {$2$};

\draw (1,4)+(0.5,0.5) node {$-\sqrt{5}$};
\draw[blue] (2,4)+(0.5,0.5) node {$\sqrt{8}$};
\draw (3,4)+(0.5,0.5) node {$1$};
\draw (4,4)+(0.5,0.5) node {$\sqrt{2}$};
\draw (5,4)+(0.5,0.5) node {$\sqrt{5}$};
\draw (1,5)+(0.5,0.5) node {$-\sqrt{8}$};
\draw (2,5)+(0.5,0.5) node {$-\sqrt{5}$};
\draw (3,5)+(0.5,0.5) node {$2$};
\draw (4,5)+(0.5,0.5) node {$\sqrt{5}$};
\draw (5,5)+(0.5,0.5) node {$\sqrt{8}$};
\draw[red,thick,{Stealth[length=3mm]}-{Stealth[length=3mm]}] (0,3.5)--(7,3.5);
\draw[red,thick,{Stealth[length=3mm]}-{Stealth[length=3mm]}] (3.5,7)--(3.5,0);
\end{tikzpicture}
\caption{$f(x,y)$}
\label{sfig:AS-F}
\end{subfigure}
\begin{subfigure}[b]{0.43\textwidth}
  \includegraphics[width=0.95\textwidth]{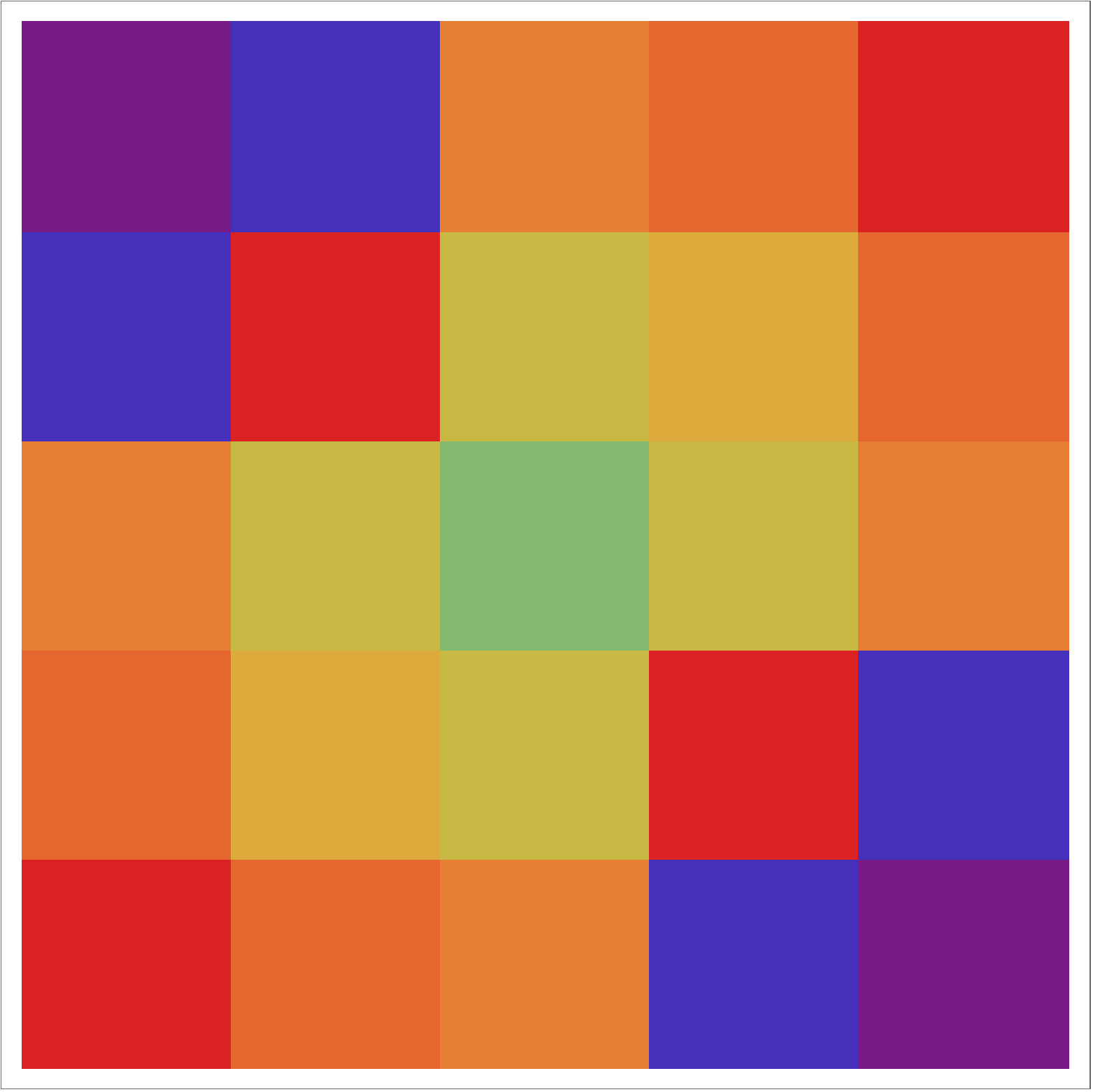}
  \vspace*{9mm}
\caption{Square in color}
\label{sfig:AS-S}
\end{subfigure}
\caption{Archeological Satorian function and square}
\label{fig:AS-FS}
\end{figure}
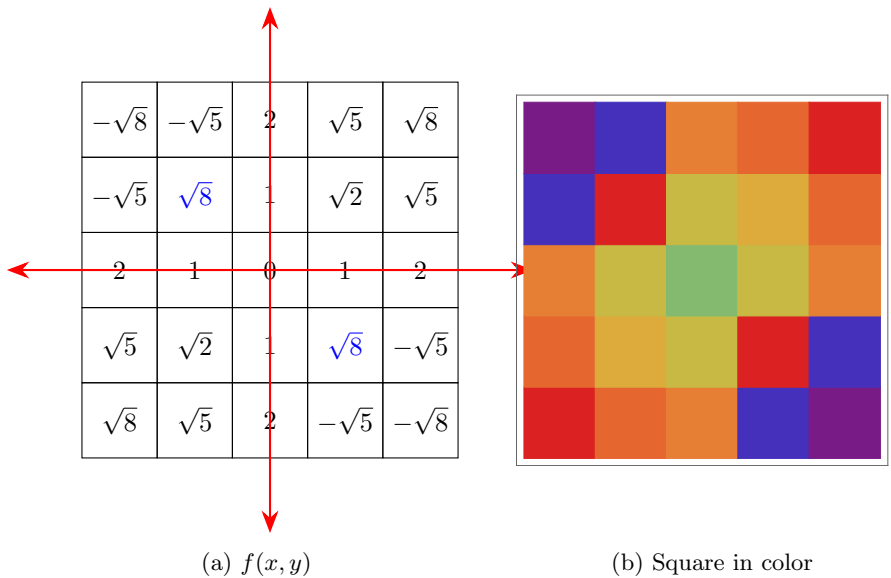

Once the definitions of a Satorian function and square have been given, it is
possible to realize that the archeological Satorian square is no more than a
particular case of a $SS_5$, where the Satorian function is given by a
discrete set.
Observe that the difference with the $SS_5$ given by Equation (\ref{eq:norm}) and
shown in Figures \ref{sfig:fxy} and \ref{sfig:SS5-C}, are the values written
in blue in Figure \ref{sfig:AS-F}.

Generating Satorian squares by using numbers is much easier than using
palyndromes or even letters.
As example, I tried to build a Sator Square of dimension $n=7$ by
using  spanish words.
To do this, it is necessary 4 words, one palindrome of 7 letters and 3 words,
which satisfy the corresponding simmetries.
Unfortunately, I found only three words\footnote{The author will appreciate if
  an example with four  words is given.}.
The palindrome  ANILINA (aniline), ECUADOR, CIANURO
(cyanide) and the no meaning word UANIMUD, see Figure \ref{fig:SS7}. 

\begin{figure}[h!]
  \centering
\begin{tikzpicture}
\foreach \x in {1,2,3,4,5,6,7}
\foreach \y in {1,2,3,4,5,6,7}
{
  \draw (\x/2,\y/2)  rectangle ++(0.5,0.5);
}
\draw (0.5,0.5)+(0.25,0.25) node {R}; 
\draw (1,0.5)+(0.25,0.25) node {O};
\draw (1.5,0.5)+(0.25,0.25) node {D};
\draw (2,0.5)+(0.25,0.25) node {\textbf{A}};
\draw (2.5,0.5)+(0.25,0.25) node {U};
\draw (3,0.5)+(0.25,0.25) node {C};
\draw (3.5,0.5)+(0.25,0.25) node {E};
\draw (0.5,1)+(0.25,0.25) node {O};
\draw (1,1)+(0.25,0.25) node {R};
\draw (1.5,1)+(0.25,0.25) node {U};
\draw (2,1)+(0.25,0.25) node {\textbf{N}};
\draw (2.5,1)+(0.25,0.25) node {A};
\draw (3,1)+(0.25,0.25) node {I};
\draw (3.5,1)+(0.25,0.25) node {C};
\draw (0.5,1.5)+(0.25,0.25) node {D};
\draw (1,1.5)  +(0.25,0.25) node {U};
\draw (1.5,1.5)+(0.25,0.25) node {M};
\draw (2,1.5)  +(0.25,0.25) node {\textbf{I}};
\draw (2.5,1.5)+(0.25,0.25) node {N};
\draw (3,1.5)  +(0.25,0.25) node {A};
\draw (3.5,1.5)+(0.25,0.25) node {U};
\draw (0.5,2)+(0.25,0.25) node {\textbf{A}};
\draw (1,2)  +(0.25,0.25) node {\textbf{N}};
\draw (1.5,2)+(0.25,0.25) node {\textbf{I}};
\draw (2,2)  +(0.25,0.25) node {\textbf{L}};
\draw (2.5,2)+(0.25,0.25) node {\textbf{I}};
\draw (3,2)  +(0.25,0.25) node {\textbf{N}};
\draw (3.5,2)+(0.25,0.25) node {\textbf{A}};
\draw (0.5,2.5)+(0.25,0.25) node {U};
\draw (1,2.5)  +(0.25,0.25) node {A};
\draw (1.5,2.5)+(0.25,0.25) node {N};
\draw (2,2.5)  +(0.25,0.25) node {\textbf{I}};
\draw (2.5,2.5)+(0.25,0.25) node {M};
\draw (3,2.5)  +(0.25,0.25) node {U};
\draw (3.5,2.5)  +(0.25,0.25) node{D};
\draw (0.5,3)+(0.25,0.25) node {C};
\draw (1,3)  +(0.25,0.25) node {I};
\draw (1.5,3)+(0.25,0.25) node {A};
\draw (2,3)  +(0.25,0.25) node {\textbf{N}};
\draw (2.5,3)+(0.25,0.25) node {U};
\draw (3,3)  +(0.25,0.25) node {R};
\draw (3.5,3)  +(0.25,0.25) node{O};
\draw (0.5,3.5)+(0.25,0.25) node {E};
\draw (1,3.5)  +(0.25,0.25) node {C};
\draw (1.5,3.5)+(0.25,0.25) node {U};
\draw (2,3.5)  +(0.25,0.25) node {\textbf{A}};
\draw (2.5,3.5)+(0.25,0.25) node {D};
\draw (3,3.5)  +(0.25,0.25) node {O};
\draw (3.5,3.5)  +(0.25,0.25) node{R};
\end{tikzpicture}
\caption{Satorian Square of dimension $n=7$. It is built with three spanish
  words, ANILINA(aniline), ECUADOR, CIANURO (cyanide) and one word without
meaning UANIMUD.}
\label{fig:SS7}
\end{figure}

\section{Satorian matrices}

By considering a Satorian Square as a matrix, it follows that it is a
bisymmetric matrix, see \cite{wang:2018}. 

\begin{definition}[Bisymmetric matrix]
An $n \times n$ matrix $A = (a_{ij} )\in \mathbb{R}$ is a bisymmetric matrix if $a_{i,j} = a_{j,i}$ and
$a_{i,j} = a_{n+1-j, n+1-i}$ for all $i,j = 1, 2, . . . , n$, that is,

\begin{equation}
\left(
\begin{array}{ccccccc}
 a_{11} & a_{12} & a_{13} & \ldots & a_{1,n-2} & a_{1,n-1} & a_{1,n} \\
 a_{12} & a_{22} & a_{23} & \ldots & a_{2,n-2} & a_{2,n-1} & a_{1,n-1} \\
 a_{13} & a_{23} & a_{33} & \ldots & a_{3,n-2} & a_{2,n-2} & a_{1,n-2} \\  
 \vdots & \vdots & \vdots & \ddots & \vdots & \vdots & \vdots \\
  a_{1,n-2} & a_{2,n-2} & a_{3,n-2} & \ldots & a_{33} & a_{23} & a_{13}\\
  a_{1,n-1} &  a_{2,n-1} & a_{2,n-2} & \ldots & a_{23} & a_{22} & a_{12}\\
a_{1,n}  & a_{1,n-1}  & a_{1,n-2} & \ldots & a_{13} & a_{12} & a_{11} \\
\end{array}
\right)
\end{equation}
\end{definition}

Bisymmetric matrices have several applications, such as numerical analysis,
graph theory or structural mechanics, vibration of beams and eigenavalue
problems see \cite{nouri:2012}, \cite{peng:2004}, \cite{peng:2007},
\cite{yuan:2013}. 

\begin{definition}[Satorian matrix]
The matrix defined by the values of a Satorian function or a Satorian
square of dimension $n$ is called a Satorian matrix and it will be noted as
$SM_n$.   
\end{definition}

A Satorian matrix $SM_n$ is a bisymmetric matrix.


\section{Satorian surfaces and solids}

If the values of a Satorian Matrix are considered as the values on the
$z$-axis, symmetric surfaces and solids can be obtained.
Thus, the Satorian matrix given by Equation (\ref{eq:SM5}) corresponds to the
Satorian Square $SS5$ shown in Figure \ref{fig:SQColor}.

\begin{equation}
  SM_5=
\left(
\begin{array}{ccccc}
 -\sqrt{8} & -\sqrt{5} & 2 & \sqrt{5} & \sqrt{8} \\
 -\sqrt{5} & -\sqrt{2} & 1 & \sqrt{2} & \sqrt{5} \\
 2 & 1 & 0 & 1 & 2 \\
 \sqrt{5} & \sqrt{2} & 1 & -\sqrt{2} & -\sqrt{5} \\
 \sqrt{8} & \sqrt{5} & 2 & -\sqrt{5} & -\sqrt{8} \\
\end{array}
\right)
\label{eq:SM5}
\end{equation}

\begin{figure}[h!]
  \centering
\begin{subfigure}[b]{0.49\textwidth}
\begin{center}
\includegraphics[width=0.95\textwidth]{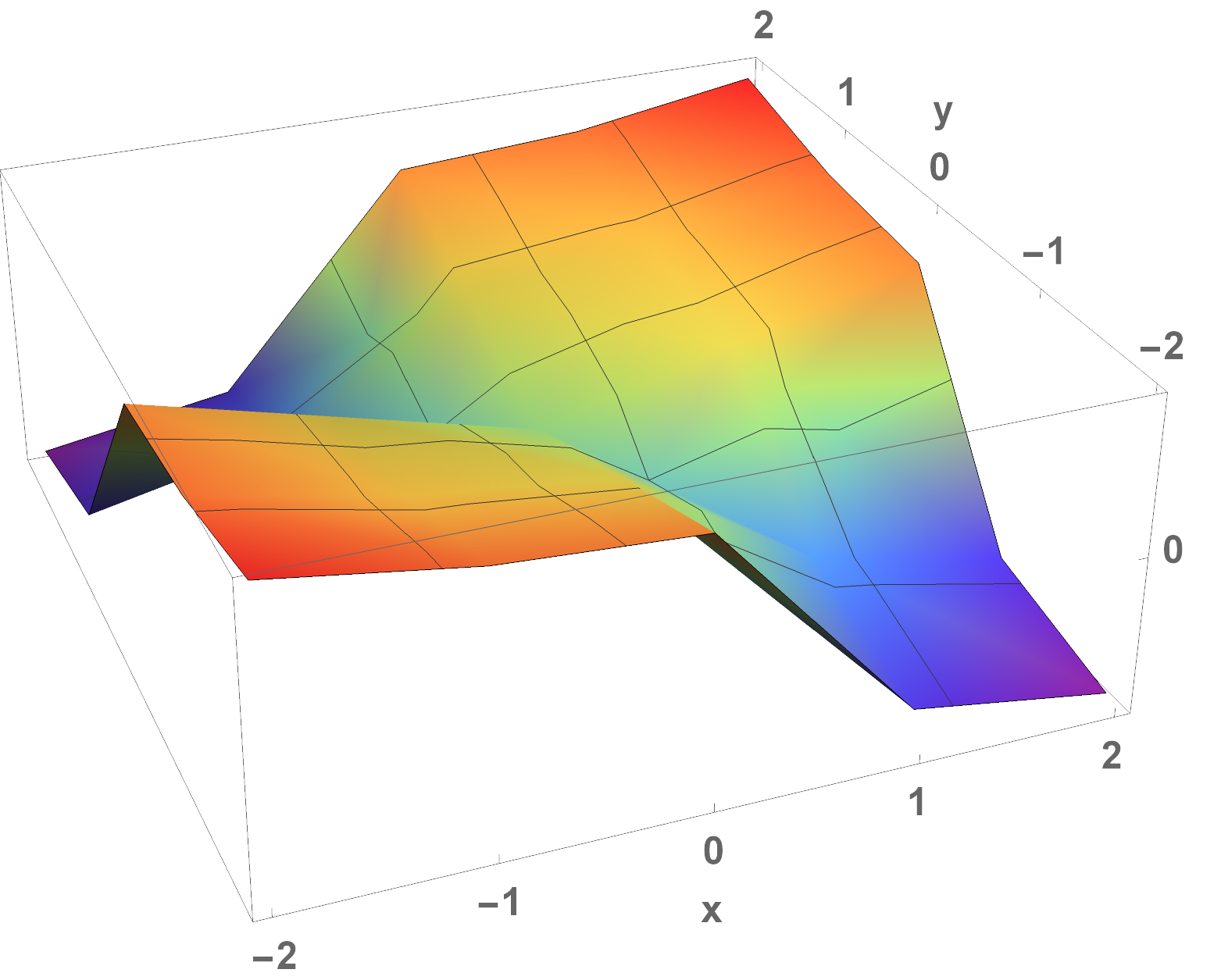}
\caption{Surface} 
\label{fig:SQ53D}
\end{center}
\end{subfigure}
\begin{subfigure}[b]{0.49\textwidth}
\begin{center}
\includegraphics[width=0.9\textwidth]{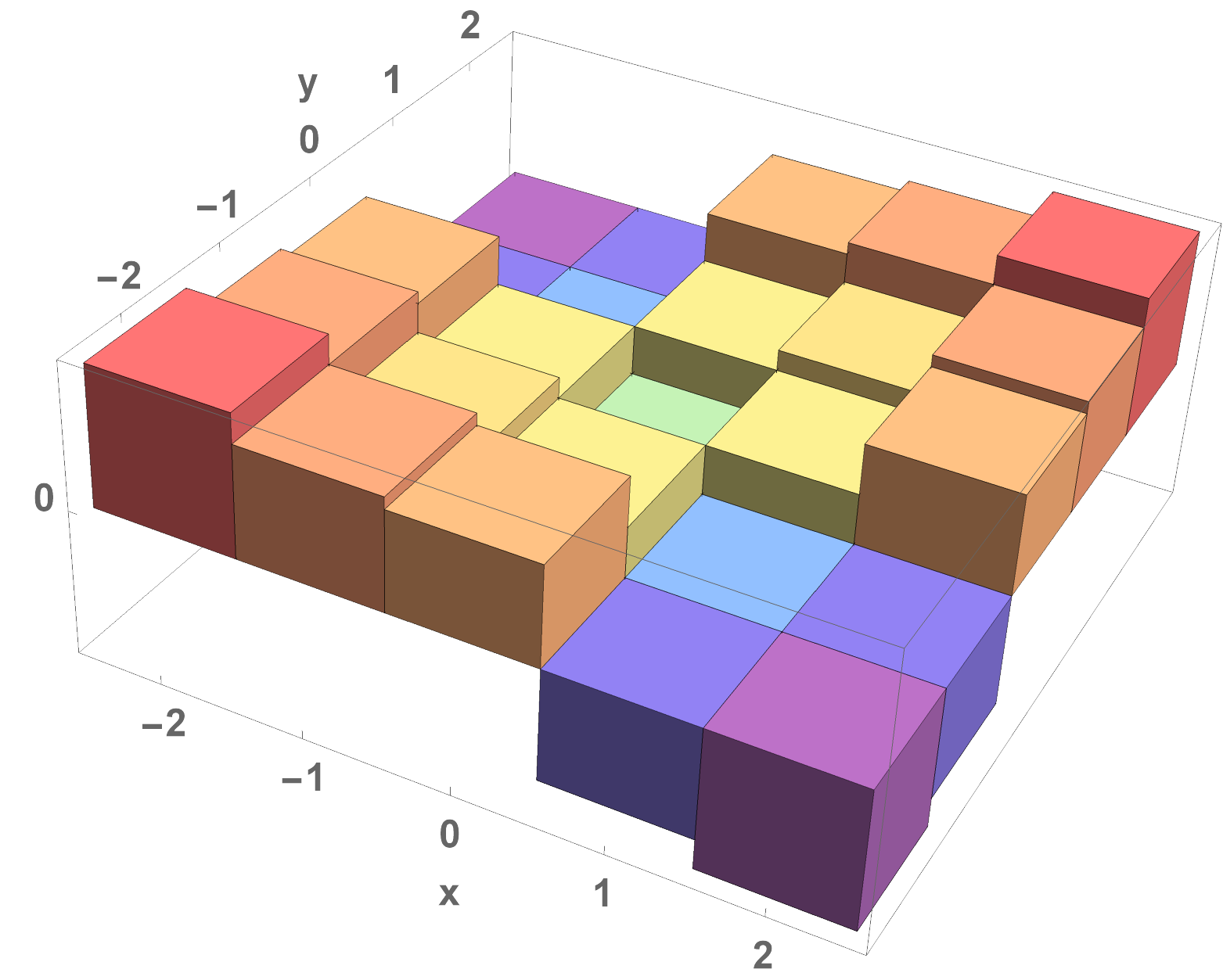}
\caption{Solid}
\label{fig:SQ5Cubes}
\end{center}
\end{subfigure}
\caption{Satorian surface and solid corresponding to the Satorian matrix
  $SM_5$}
\label{fig:SS3D}
\end{figure}

Observe that the matrix $SM_5$ depends on the function $f(x,y)$ given by
Equation (\ref{eq:norm}) and shown in Figure \ref{sfig:fxy}.
The Figure \ref{fig:SS3D} shows the surface and solid corresponding to the
Matrix $SM_5$.

\begin{figure}[h!]
  \centering
\begin{subfigure}[b]{0.49\textwidth}
\begin{center}
\includegraphics[width=0.95\textwidth]{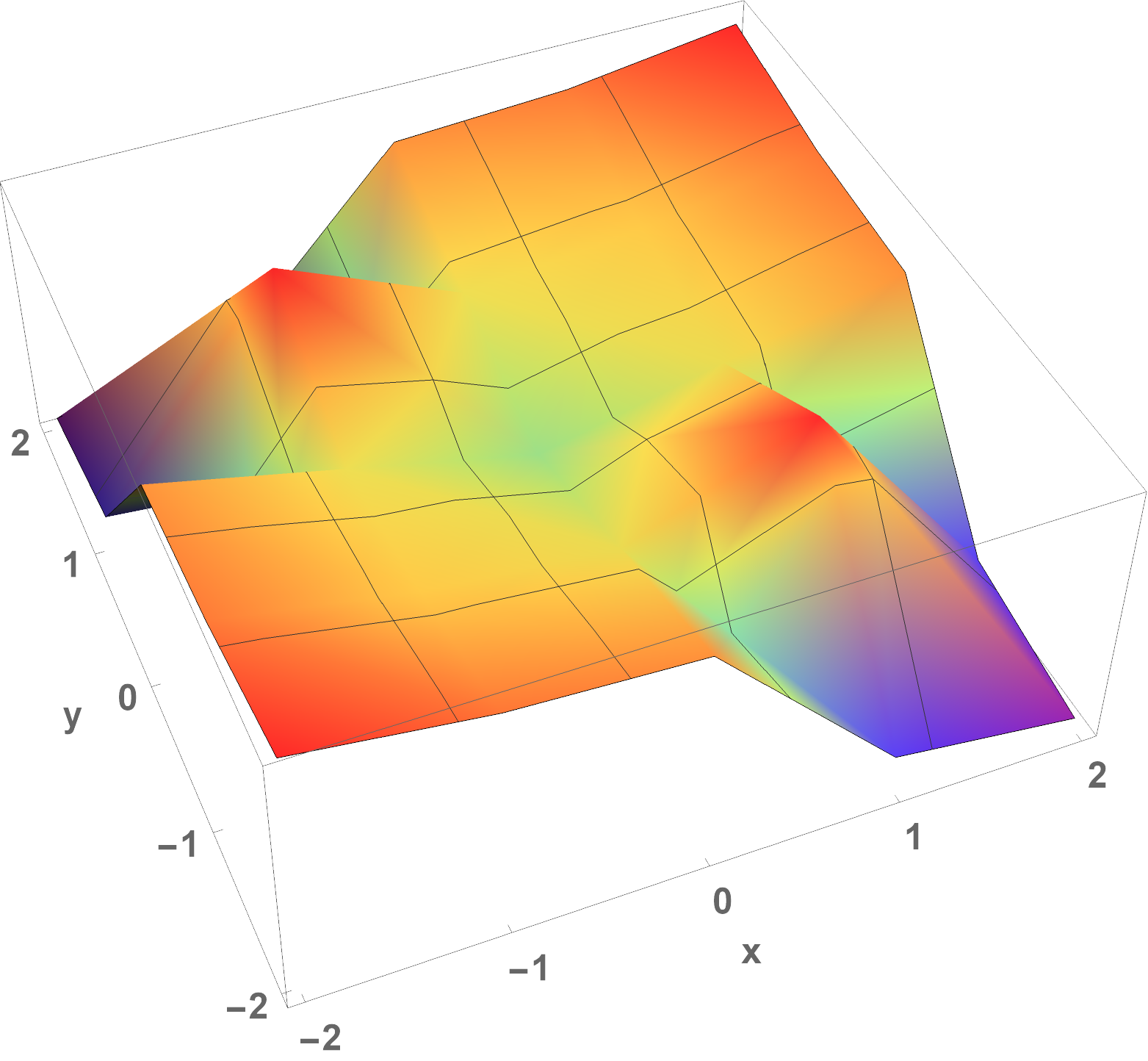}
\caption{Surface} 
\label{fig:SQ53D}
\end{center}
\end{subfigure}
\begin{subfigure}[b]{0.49\textwidth}
\begin{center}
\includegraphics[width=0.9\textwidth]{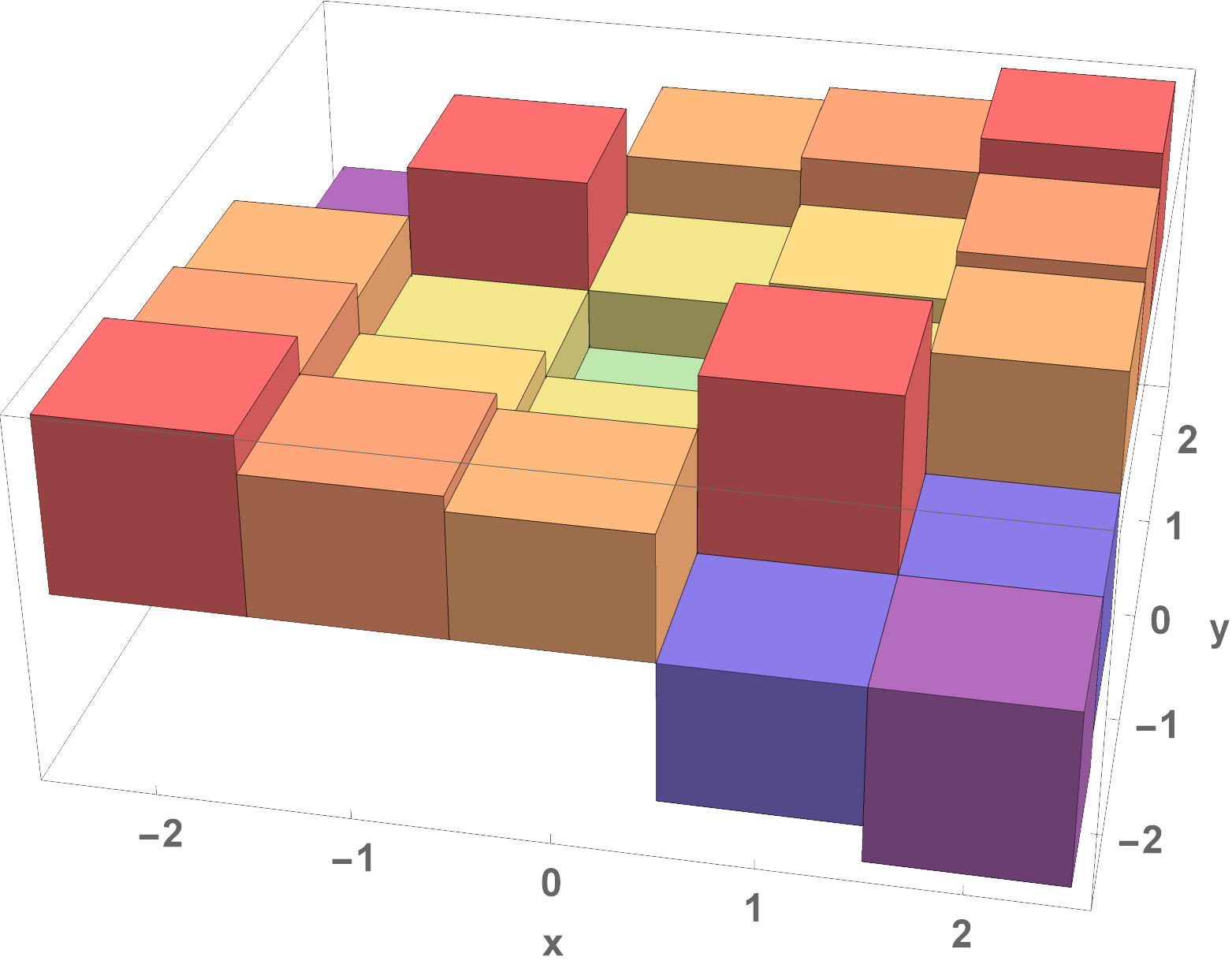}
\caption{Solid}
\label{fig:SQ5Cubes}
\end{center}
\end{subfigure}
\caption{Satorian surface and solid corresponding to the
  Archeological Satorian matrix $ASM_5$}
\label{fig:ASS3D}
\end{figure}

The matrix corresponding to the Archeological Satorian Square is given in
Equation (\ref{eq:ASM5}), and the surface and solid corresponding to the
Archeological Satorian square are shown in Figure \ref{fig:ASS3D}

\begin{equation}
  ASM_5=
\left(
\begin{array}{ccccc}
 -\sqrt{8} & -\sqrt{5} & 2 & \sqrt{5} & \sqrt{8} \\
 -\sqrt{5} & \sqrt{8} & 1 & \sqrt{2} & \sqrt{5} \\
 2 & 1 & 0 & 1 & 2 \\
 \sqrt{5} & \sqrt{2} & 1 & \sqrt{8} & -\sqrt{5} \\
 \sqrt{8} & \sqrt{5} & 2 & -\sqrt{5} & -\sqrt{8} \\
\end{array}
\right)
\label{eq:ASM5}
\end{equation}

As it was mentioned before, the Satorian functions can be defined by a
discrete set of points or by a set of real functions.
Thus, let us consider an application of a Gaussian function given by

\begin{equation}
  f(x,y)=\left\{
    \begin{array}{rl}
      |x+y| & \text{if } xy=0,\\\\
      \exp{(-x^2-y^2)} & \text{if } xy>0,\\\\
      -\exp{(-x^2-y^2)} & \text{if } xy<0
    \end{array} \right.
  \label{eq:gauss}
 \end{equation}

\begin{figure}[H]
\centering
\begin{subfigure}[b]{0.44\textwidth}
\begin{tikzpicture}
\foreach \x in {1,2,3,4,5}
\foreach \y in {1,...,5}
{
  \draw (\x,\y)  rectangle ++(1,1);
}
\draw (1,1)+(0.5,0.5) node {$\cdot$}; 
\draw (2,1)+(0.5,0.5) node {$\cdot$};
\draw (3,1)+(0.5,0.5) node {-1};
\draw (4,1)+(0.5,0.5) node {$\cdot$};
\draw (5,1)+(0.5,0.5) node {$\cdot$};
\draw (1,2)+(0.5,0.5) node {$\cdot$};
\draw (2,2)+(0.5,0.5) node {$\cdot$};
\draw (3,2)+(0.5,0.5) node {-0.5};
\draw (4,2)+(0.5,0.5) node {$\cdot$};
\draw (5,2)+(0.5,0.5) node {$\cdot$};
\draw (1,3)+(0.5,0.5) node {-1};
\draw (2,3)+(0.5,0.5) node {-0.5};
\draw (3,3)+(0.5,0.5) node {0};
\draw (4,3)+(0.5,0.5) node {0.5};
\draw (5,3)+(0.5,0.5) node {1};
\draw (1,4)+(0.5,0.5) node {$\cdot$};
\draw (2,4)+(0.5,0.5) node {$\cdot$};
\draw (3,4)+(0.5,0.5) node {0.5};
\draw (4,4)+(0.5,0.5) node {$\cdot$};
\draw (5,4)+(0.5,0.5) node {$\cdot$};
\draw (1,5)+(0.5,0.5) node {$\cdot$};
\draw (2,5)+(0.5,0.5) node {$\cdot$};
\draw (3,5)+(0.5,0.5) node {1};
\draw (4,5)+(0.5,0.5) node {$\cdot$};
\draw (5,5)+(0.5,0.5) node {$\cdot$};
\draw[red,{Stealth[length=3mm]}-{Stealth[length=3mm]}] (0,3.5)--(7,3.5);
\draw[red,{Stealth[length=3mm]}-{Stealth[length=3mm]}] (3.5,7)--(3.5,0);
\draw[red] (6.8,3.2) node {$x$};
\draw[red] (3.8,6.8) node {$y$};
\end{tikzpicture}
\caption{$(x,y)$}
\label{sfig:xy2}
\end{subfigure}
\hfill
\begin{subfigure}[b]{0.44\textwidth}
\begin{tikzpicture}
\foreach \x in {1,2,3,4,5}
\foreach \y in {1,...,5}
{
  \draw (\x,\y)  rectangle ++(1,1);
}
\draw (1,1)+(0.5,0.5) node {{\small 0.135}}; 
\draw (2,1)+(0.5,0.5) node {{\small0.287}};
\draw (3,1)+(0.5,0.5) node {{\small 1}};
\draw (4,1)+(0.5,0.5) node {{\small-0.287}};
\draw (5,1)+(0.5,0.5) node {{\small-0.135}};
\draw (1,2)+(0.5,0.5) node {{\small 0.287}};
\draw (2,2)+(0.5,0.5) node {{\small 0.607}};
\draw (3,2)+(0.5,0.5) node {{\small 0.5}};
\draw (4,2)+(0.5,0.5) node {{\small-0.607}};
\draw (5,2)+(0.5,0.5) node {{\small -0.287}};
\draw (1,3)+(0.5,0.5) node {{\small 1}};
\draw (2,3)+(0.5,0.5) node {{\small 0.5}};
\draw (3,3)+(0.5,0.5) node {{\small 0}};
\draw (4,3)+(0.5,0.5) node {{\small 0.5}};
\draw (5,3)+(0.5,0.5) node {{\small 1}};
\draw (1,4)+(0.5,0.5) node {{\small -0.287}};
\draw (2,4)+(0.5,0.5) node {{\small -0.607}};
\draw (3,4)+(0.5,0.5) node {{\small 0.5}};
\draw (4,4)+(0.5,0.5) node {{\small 0.607}};
\draw (5,4)+(0.5,0.5) node {{\small 0.287}};
\draw (1,5)+(0.5,0.5) node {{\small -0.135}};
\draw (2,5)+(0.5,0.5) node {{\small -0.287}};
\draw (3,5)+(0.5,0.5) node {{\small 1}};
\draw (4,5)+(0.5,0.5) node {{\small 0.287}};
\draw (5,5)+(0.5,0.5) node {{\small 0.135}};
\draw[red,{Stealth[length=3mm]}-{Stealth[length=3mm]}] (0,3.5)--(7,3.5);
\draw[red,{Stealth[length=3mm]}-{Stealth[length=3mm]}] (3.5,7)--(3.5,0);
\draw[red] (6.8,3.2) node {$x$};
\draw[red] (3.8,6.8) node {$y$};
\end{tikzpicture}
\caption{$f(x,y)$}
\label{sfig:gxy}
\end{subfigure}
\caption{Application of the Equation (\ref{eq:gauss}) on a discrete square}
\label{fig:gauss-n5}
\end{figure}
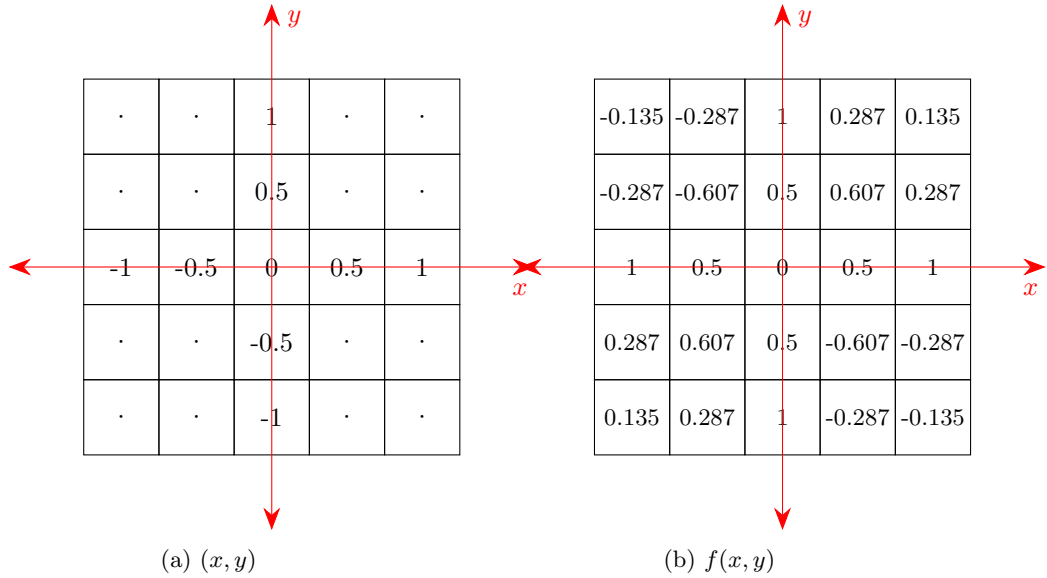

 By aplying the Satorian function given in Equation (\ref{eq:gauss}) in a square
 of dimension $n=5$ the corresponding matrix is obtained, see Figure
 \ref{fig:gauss-n5}.

By assigning colours to the Satorian square shown in Figure \ref{sfig:gxy} and
by plotting the values on the $z$ axis, the corresponding square and surface
are obtained. 

\begin{figure}[H]
  \centering
\begin{subfigure}[b]{0.49\textwidth}
\begin{center}
\includegraphics[width=0.7\textwidth]{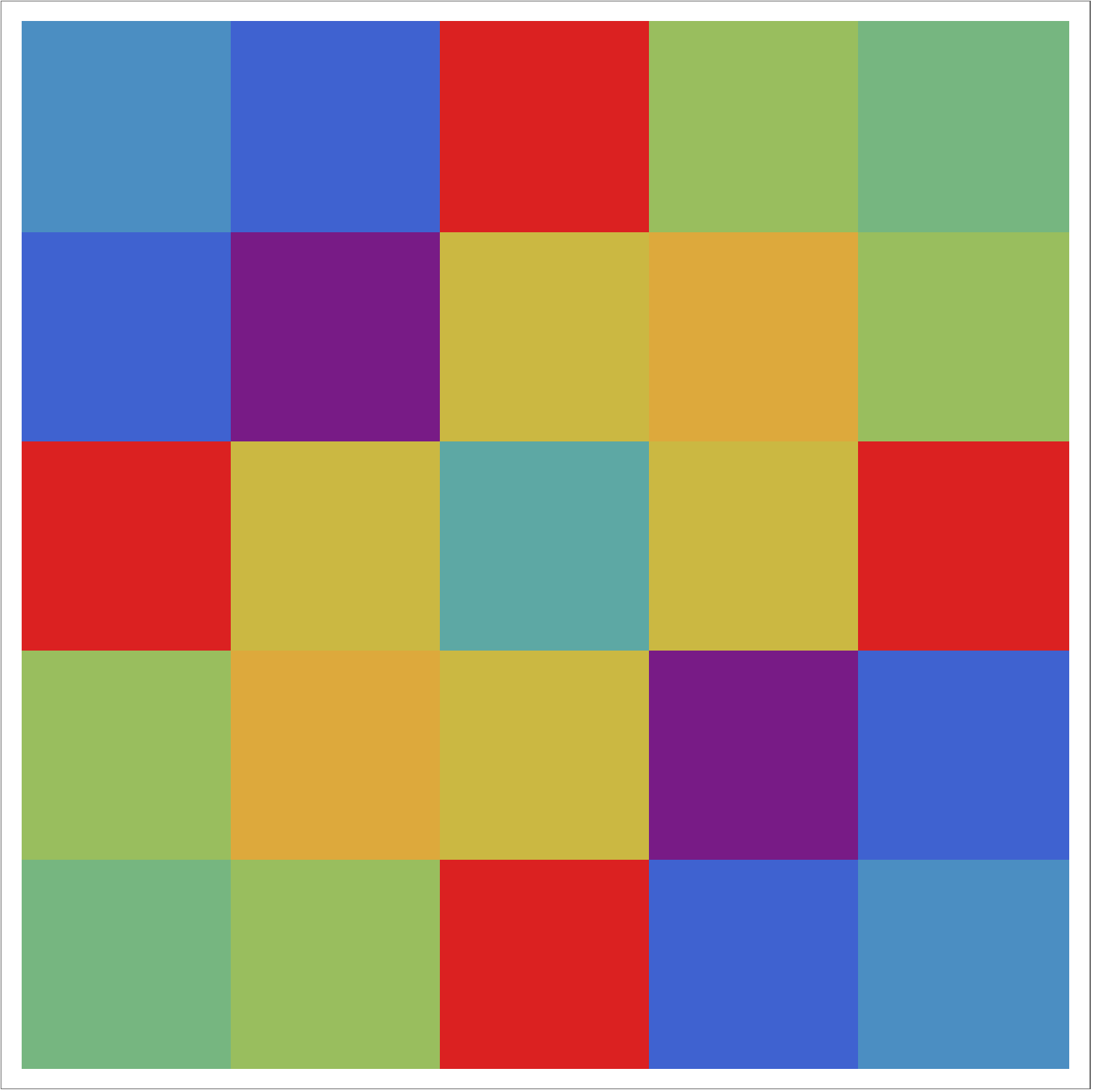}
\caption{Square}
\label{fig:GS}
\end{center}
\end{subfigure}
\begin{subfigure}[b]{0.49\textwidth}
\begin{center}
\includegraphics[width=0.95\textwidth]{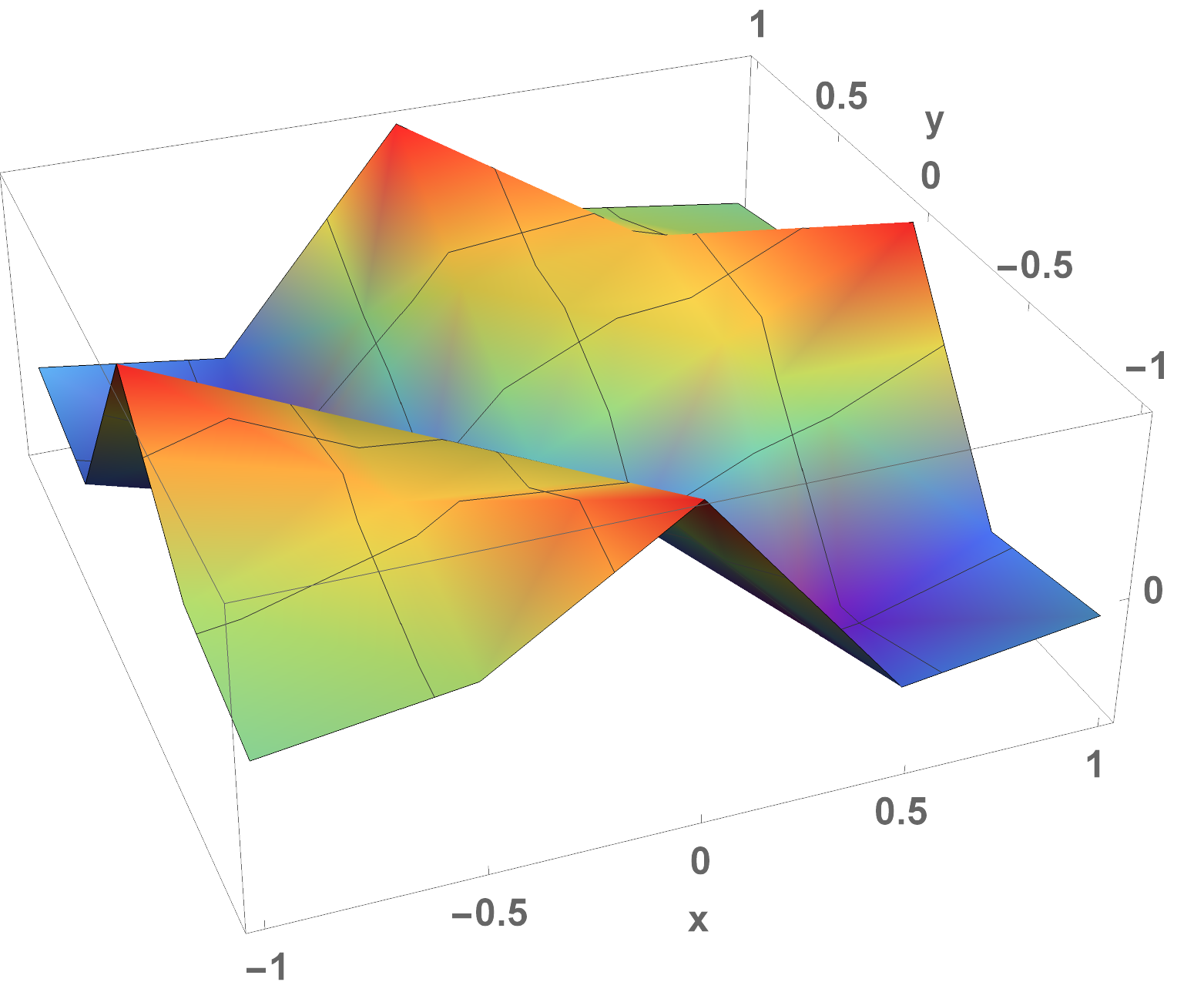}
\caption{Surface}
\label{fig:SQ5Cubes}
\end{center}
\end{subfigure}
\caption{Satorian square and surface corresponding to the function
  given by Equation (\ref{eq:gauss})}
\label{fig:GS3D}
\end{figure}

\begin{figure}[H]
  \centering
\begin{subfigure}[b]{0.49\textwidth}
\begin{center}
\includegraphics[width=0.85\textwidth]{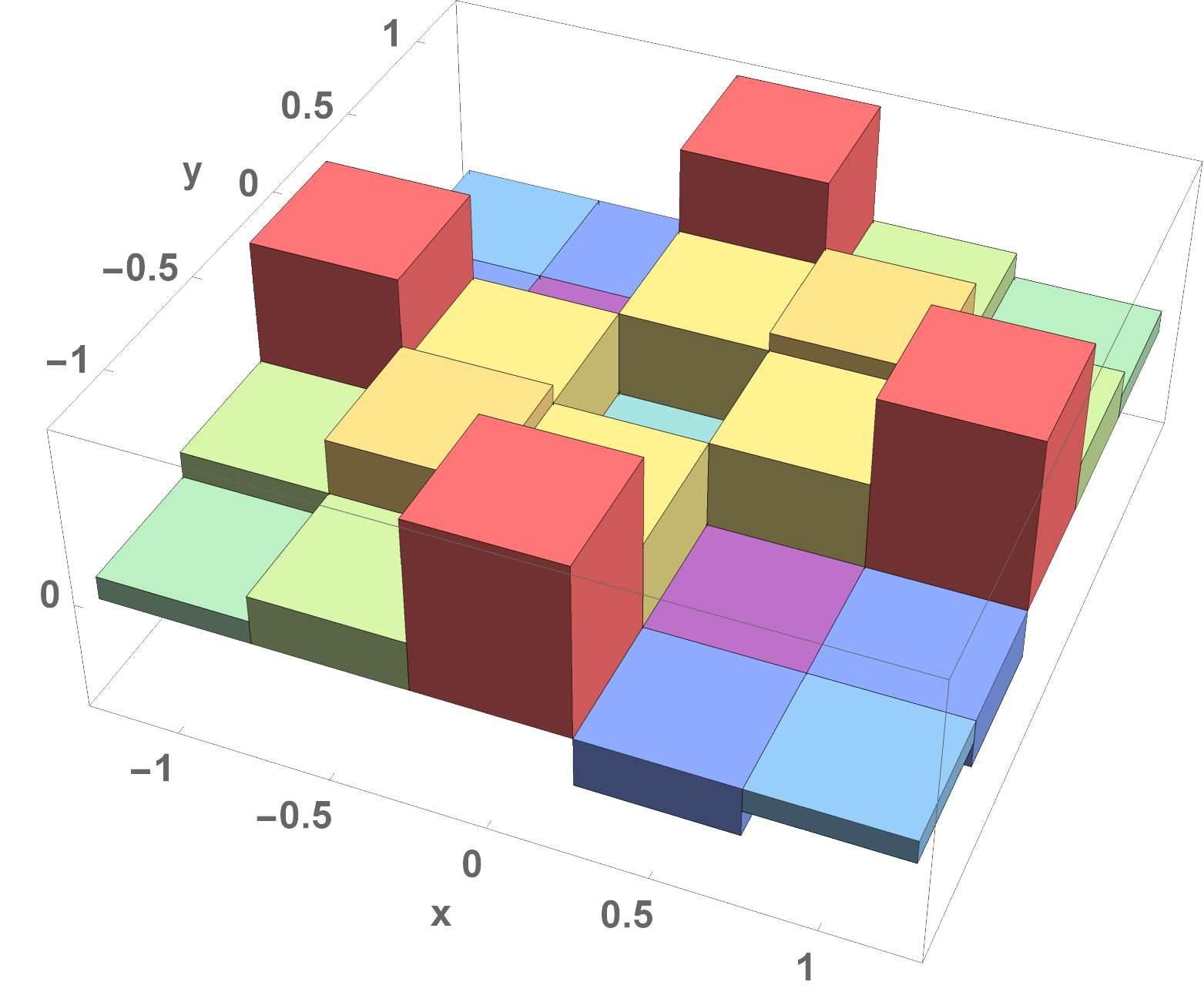}
\caption{Solid}
\label{fig:GSolid}
\end{center}
\end{subfigure}
\begin{subfigure}[b]{0.49\textwidth}
\begin{center}
\includegraphics[width=0.95\textwidth]{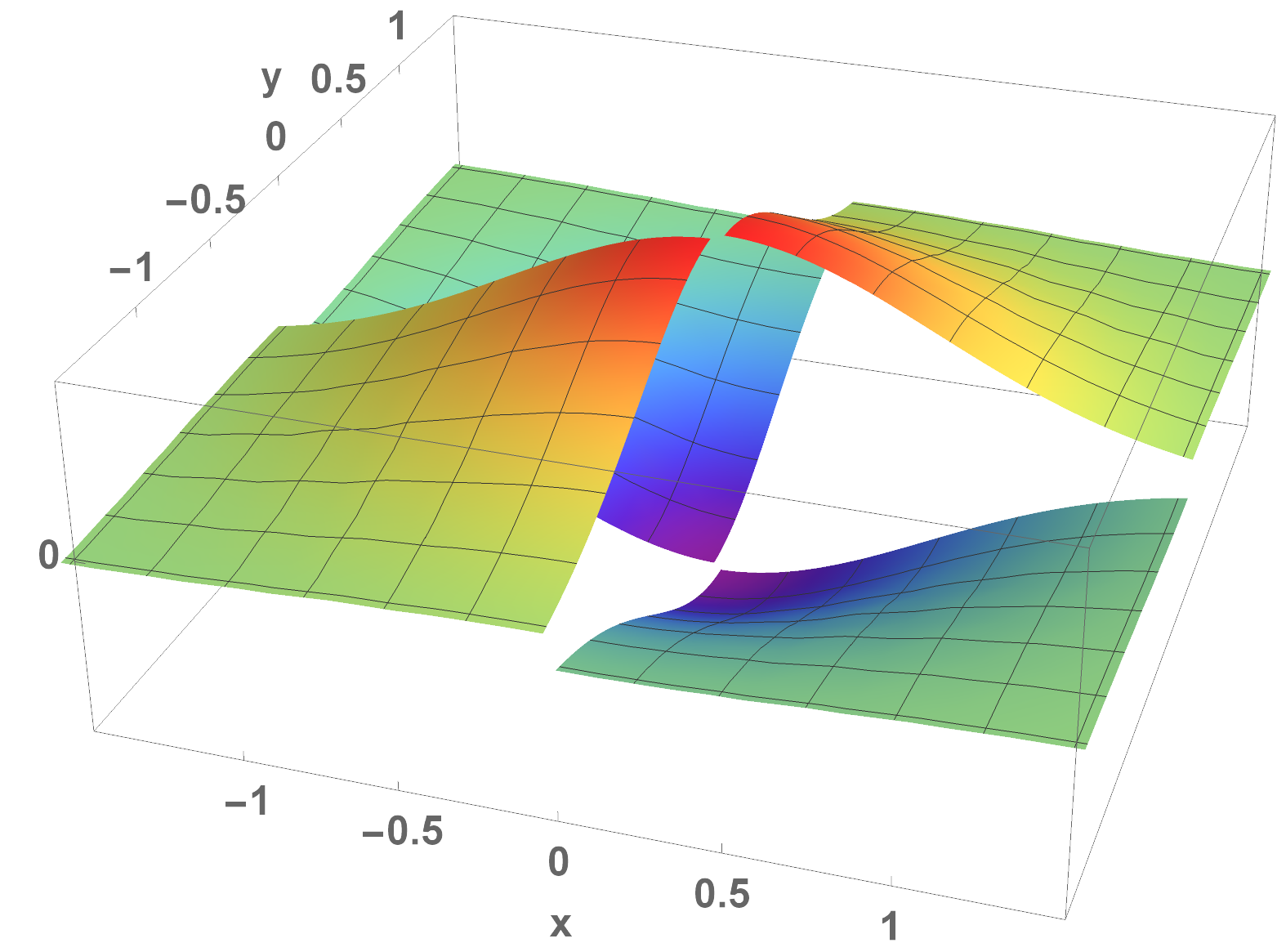}
\caption{Surface in $\mathbb{R}$}
\label{fig:SSUrfaceR}
\end{center}
\end{subfigure}
\caption{Satorian solid and surface corresponding to the function
  given by Equation (\ref{eq:gauss})}
\label{fig:GS3D}
\end{figure}

The Figure \ref{fig:GSolid} shows the solid correspondig to the  function
given by Equation (\ref{eq:gauss}) and the Figure \ref {fig:SSUrfaceR} shows
the surface when the function $f(xy)$ is evaluated in  $x,y\in \mathbb{R}$.

\section{Arabic architecture. The Alhambra}
\begin{figure}[h!]
\centering
\includegraphics[width=0.95\textwidth]{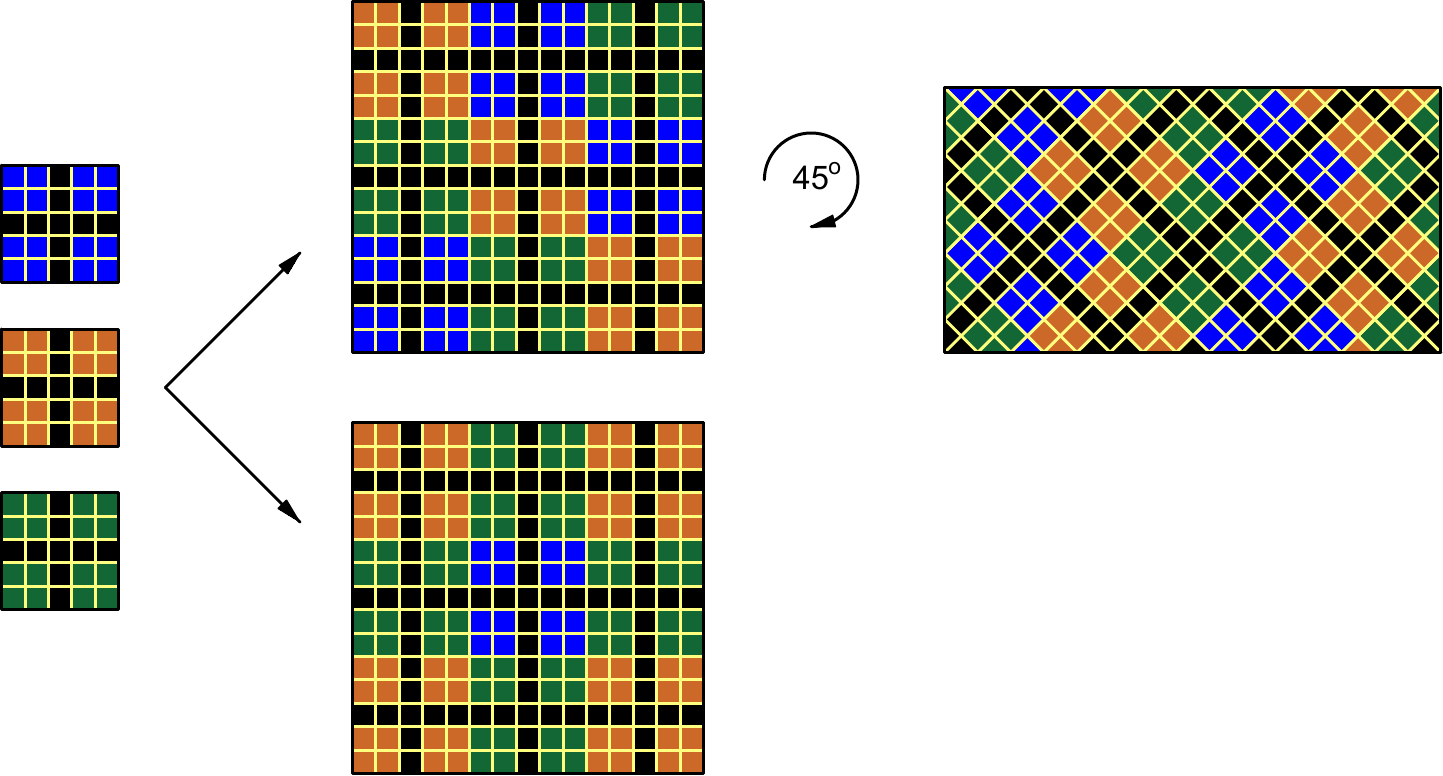}
\caption{Creation of a Satorian mosaic} 
\label{fig:Al_01}
\end{figure}

\begin{figure}[H]
\centering
\includegraphics[width=0.5\textwidth]{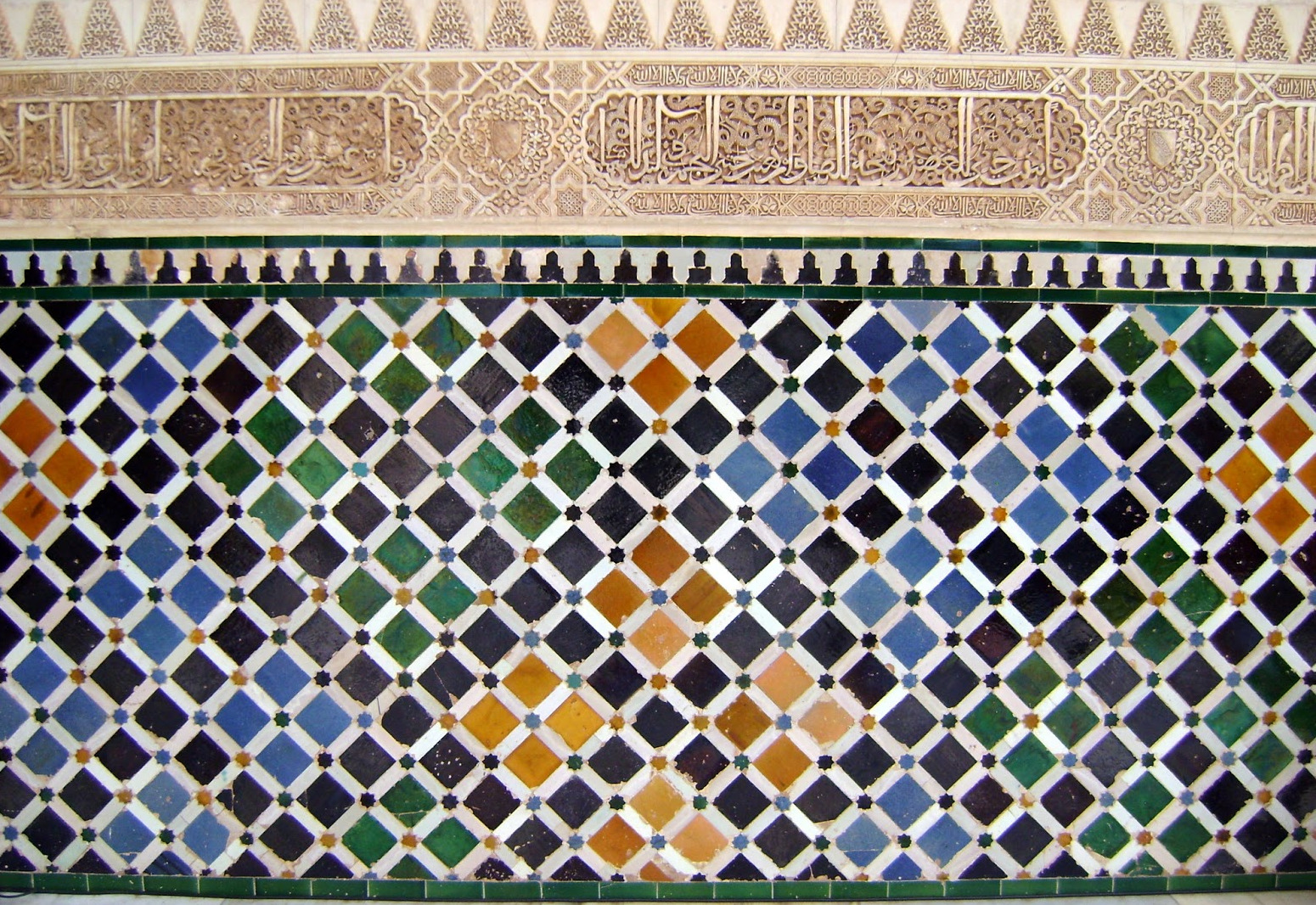}
\caption{Alicatado on one wall at the Court of the Myrtles from the Comares
  Palace inside the Alhambra}  
\label{fig:Al_02}
\end{figure}

The kind of symmetry observed in the Sator Square has been also found in the
alicatados\footnote{Mosaic formed of polygonal, coloured glazed
    tiles. Made up into geometric patterns.} used for paving Spanish and
  Moorish patios but also for wall surfaces. 
The Figure \ref{fig:Al_01} show how a Satorian mosaic can be creted from three
Satorian squares $SS_5$. 
To do this, nine squares $SS_5$ are joined to form a bigger square, which does
not satisfy the simmetry $S_2$ and afterwards, the design is turned $45^o$.
However, by combining properly nine squares $SS_{5}$, a Satorian square
$SS_{15}$ can be obtained as well.\\
The real picture of the Satorian mosaic in the Alhambra is showed in Figure
\ref{fig:Al_02}.

A fascinating work about simmetry in the Alhambra and other fields can be found in
\cite{satoy:2008}.


\section{Conclusions}

For centuries the Sator Square has been considered as a mysterious riddle.
Most of the studies about this archeological object have been based on
religions or pseudo-sicentific disciplines, such as Christianity, Judaism,
Satanism, Philology, esotericism, magic and mysticism.
Nobody has been able to propose a rational and convincing explanation about
its meaning.\\ 

In this work, the author proposed to see this square as a simple symmetric
game of words.
However, as every symmetric object, its structure also can be related to
mathematical objects.
In fact, the structure of the Sator Square can be described as a function, a
matrix, a surface or a solid.
By defining the mathematical properties of a Sator Square, it is possible to
create Satorian structures of higher dimension.\\

Perhaps if the Sator square is analized or studied from an analytical
perspective, a reliable explanation about its origin and meaning can be
acquiered.


\bibliographystyle{ieeetr}
\bibliography{Sator_references}%

\end{document}